\documentclass[12pt]{elsarticle}

\usepackage{graphics,graphicx,subfigure,amssymb}

\def\xsum{\mathop{\sum\nolimits'}}

\journal{Applied Mathematics and Computation}

\begin{document}

\begin{frontmatter}

\title{Generation of off-critical zeros for hypercubic Epstein zeta-functions}

\author{Igor Trav\v{e}nec\corref{mycorrespondingauthor}} 
\cortext[mycorrespondingauthor]{Corresponding author}
\ead{fyzitrav@savba.sk}

\author{Ladislav \v{S}amaj}
\address{Institute of Physics, Slovak Academy of Sciences, 
D\'ubravsk\'a cesta 9, 84511 Bratislava, Slovakia}

\begin{abstract}
We study the Epstein zeta-function formulated on the $d$-dimensional
hypercubic lattice $\zeta^{(d)}(s) = \frac{1}{2} \xsum_{n_1,\ldots,n_d} 
(n_1^2+\cdots+n_d^2)^{-s/2}$ where the real part $\Re(s)>d$
and the summation runs over all integers except of the origin $(0,0,\ldots,0)$. 
An analytical continuation of the Epstein zeta-function to the whole complex
$s$-plane is constructed for the spatial dimension $d$ being a continuous 
variable ranging from $0$ to $\infty$.
Zeros of the Epstein zeta-function $\rho=\rho_x+{\rm i}\rho_y$ are defined by 
$\zeta^{(d)}(\rho) = 0$.
The nontrivial zeros split into the ``critical'' zeros (on the critical line) 
with $\rho_x=\frac{d}{2}$ and the ``off-critical'' zeros 
(off the critical line) with $\rho_x \ne \frac{d}{2}$.
Numerical calculations reveal that the critical zeros form closed or semi-open 
curves $\rho_y(d)$ which enclose disjunctive regions of the plane
$\left( \rho_x= \frac{d}{2},\rho_y \right)$.
Each curve involves a number of left/right edge points 
$\rho^*= \left( \frac{d^*}{2},\rho_y^* \right)$, 
defined by a divergent tangent ${\rm d}\rho_y/{\rm d}d\vert_{\rho^*}$.
Every edge point gives rise to two conjugate tails of off-critical zeros with 
continuously varying dimension $d$ which exhibit a singular expansion 
around the edge point, in analogy with critical phenomena for second-order 
phase transitions.
For each dimension $d>9.24555\ldots$ there exists a conjugate 
pair of {\em real} off-critical zeros which tend to the boundaries 
$0$ and $d$ of the critical strip in the limit $d\to\infty$. 
As a by-product of the formalism, we derive an exact formula for
$\lim_{d\to 0} \zeta^{(d)}(s)/d$.
An equidistant distribution of critical zeros along the imaginary axis
is obtained for large $d$, with spacing between the nearest-neighbour
zeros vanishing as $2\pi/\ln d$ in the limit $d\to\infty$.
\end{abstract}

\begin{keyword}
Epstein zeta-function; hypercubic lattice; analytic continuation; 
Jacobi theta function; zeros; critical phenomena
\MSC[2010]{11E45}
\end{keyword}

\end{frontmatter}


\section{Introduction} \label{Sec1}
Let two particles at distance $r$ interact via the Riesz potential 
$1/r^s$ with real $s$ \cite{Brau}. 
If the particles are placed equidistantly (with unit lattice spacing) on 
an infinite line and interact pairwisely by the Riesz potential, 
the energy per particle is given by the Riemann zeta-function
\cite{Riemann1859}
\begin{equation}
\zeta(s) = \frac{1}{2} \xsum_{n=-\infty}^{\infty} 
\frac{1}{\vert n\vert^s} = \sum_{n=1}^{\infty} \frac{1}{n^s}  \qquad s>1 ,
\end{equation}
where the prefactor $\frac{1}{2}$ comes from the fact that each interaction
energy is shared by a pair of particles and the prime in the first sum means 
omission of the self-energy $n=0$ term from the set of integers $n$.
The function $\zeta(s)$ can be analytically continued to the punctured plane
$\mathbb{C} \setminus \{ 1\}$.
It has a simple pole at $s=1$.
The Riemann zeta-function plays a fundamental role in the algebraic and
analytic number theories
\cite{Hardy14,Riesz16,Hardy21,Hutchinson25,Titchmarsh35,Selberg46},
see monographs \cite{Edwards74,Ivic85,Titchmarsh88}.
The Riemann hypothesis about the location of its nontrivial zeros exclusively 
on the critical line $\Re(s)=\frac{1}{2}$ (the symbol $\Re$ means 
the real part) is one of the Hilbert and Clay Millennium Prize problems.
Throughout the present paper we assume that the Riemann hypothesis holds.
The zeros of the Riemann zeta-function are tabulated in the symbolic
language \textit{Mathematica} under the symbol ZetaZero[$n$] where
the positive integer $n=1,\ldots,10^7$ denotes the $n^{\rm th}$ zero 
in the first quadrant and the negative integer $n=-1,\ldots,-10^7$ 
corresponds to its complex conjugate.
The zeta-function and its Hurwitz, Barnes \cite{Barnes04}, Epstein
\cite{Epstein03,Epstein07}, etc. generalisations have numerous applications
in mathematics (prime numbers, applied statistics 
\cite{Apostol76,Apostol90,Borwein13}) and physics (dynamical systems 
\cite{Ruelle02}, regularization in Quantum Field Theory
\cite{Elizalde94}, Casimir effect \cite{Milton01,Edery06}, 
Bose-Einstein condensation \cite{Dalvovo99}, 
see also books \cite{Kirsten10,Elizalde12}).

If the particles sit on the vertices of the $d$-dimensional hypercubic
lattice with unit spacing, the energy per particle is given by 
the hypercubic Epstein zeta-function \cite{Epstein03,Epstein07,Chowla49}
\begin{equation} \label{zetad}
\zeta^{(d)}(s) = \frac{1}{2} \xsum_{n_1,n_2,\ldots, n_d}
\frac{1}{(n_1^2+n_2^2+\ldots+n_d^2)^{s/2}}  \qquad \Re(s) > d ,
\end{equation} 
where the self-energy term $(0,\ldots,0)$ is excluded from the summation
over all integers and the spatial dimension $d$ is a positive integer.
This function can be analytically continued (regularized) to the critical 
strip $0<\Re(s)<d$ by various methods
\cite{Epstein03,Epstein07,Ennola64,Elizalde89a,Blanc15}.
One of the methods is based on the fact that if the particles charges change 
signs periodically so that the system is electrically neutral, the appropriate 
lattice sums converge for all $\Re(s)>0$ except for $s=d$ \cite{Borwein13}.
The critical line is defined by $\Re(s)=\frac{d}{2}$.
Since the hypercubic lattice is self-dual, one can use a functional relation 
that connects the lattice sums for $s\leftrightarrow d-s$ and thus get 
an analytic continuation from $\Re(s)>d$ to $\Re(s)<0$ \cite{Borwein13}, 
although the region $\Re(s)<0$ is problematic from the point of view of physical
applications. 

We are interested in zeros of the Epstein zeta-function
$\rho=\rho_x+{\rm i}\rho_y$ defined by $\zeta^{(d)}(\rho) = 0$.
Besides the trivial zeros at $\rho=-2,-4,\ldots$, there exist non-trivial
zeros which split into two sets: the ``critical'' zeros (on the critical line) 
with $\rho_x = \frac{d}{2}$ and the ``off-critical'' zeros (off the critical 
line) with $\rho_x \ne \frac{d}{2}$.

For dimensions $d=2$, $4$, $6$ and $8$, the hypercubic Epstein zeta-functions 
can be expressed in terms of certain one-dimensional sums 
\cite{Glasser73a,Glasser73b,Zucker74,Borwein14}.
In the general analysis of these dimensions the Riemann hypothesis
will be assumed to be valid.

For the square lattice ($d=2$), it holds that \cite{Lorenz71,Hardy19} 
\begin{equation} \label{zeta2}
\zeta^{(2)}(s) = 2\ \zeta\left(\frac{s}{2}\right)
\beta\left(\frac{s}{2}\right),
\end{equation}
where
\begin{equation} 
\beta\left(\frac{s}{2}\right) = \sum_{k=0}^\infty\frac{(-1)^k}{(2k+1)^{s/2}}
= \frac{1}{2^s}\left[\zeta\left(\frac{s}{2},\frac{1}{4}\right)
-\zeta\left(\frac{s}{2},\frac{3}{4}\right)\right] 
\end{equation}
is the Dirichlet beta-function which is a special case of Dirichlet $L$-series
\cite{Borwein13}; here, $\zeta(s/2,a)=\sum_{n=0}^\infty(n+a)^{-s/2}$
denotes the Hurwitz zeta-function.
Provided that the Riemann hypothesis holds for the zeta function, zeros of
$\beta(s/2)$ are localized on the critical line $\Re(s/2)=\frac{1}{2}$
\cite{Lander18}, so that all nontrivial zeros of $\zeta^{(2)}(s)$ 
are constrained to the critical line $\Re(s)=1$.
The statistics of gaps between zeros was studied numerically \cite{Hejhal87}
as well as analytically \cite{Bogomolny94}.
The anisotropic (rectangular) lattice sums are of special interest
due to the presence of zeros off the critical line
\cite{Potter35,Davenport36,Bateman64,Stark67,McPhedran16}.
The distribution of critical zeros for general two-dimensional periodic 
structures was studied in \cite{Jutila05,Baier17}. 

For $d=4$, the Epstein zeta-function is expressible as
\cite{Zucker74,Borwein14}
\begin{equation} \label{zeta4}
\zeta^{(4)}(s) =  4 \left(1-2^{2-s}\right) \zeta\left(\frac{s}{2}\right)
\zeta\left(\frac{s}{2}-1\right) .
\end{equation}
The critical zeros are given by $2^{2-\rho}=1$, i.e. 
$\rho = 2 + 2\pi{\rm i}k/\ln 2$ $(k=0,\pm 1,\pm 2,\ldots)$.
There are also off-critical zeros lying on the lines $\rho_x=1$ and $\rho_x=3$.

For $d=6$, it holds that \cite{Zucker74,Borwein14}
\begin{equation} \label{zeta6}
\zeta^{(6)}(s)= 8\,  \beta\left(\frac{s}{2}\right)
\zeta\left(\frac{s}{2}-2\right)-2\, \beta\left(\frac{s}{2}-2\right)
\zeta\left(\frac{s}{2}\right).
\end{equation}
Besides the critical zeros on the axis $\Re(s)=3$, there exist also 
off-critical zeros dispersed in the complex plane.

For $d=8$ \cite{Zucker74,Borwein14}, 
\begin{equation} \label{zeta8}
\zeta^{(8)}(s)= 8 \left(1-2^{1-s/2}+4^{2-s/2}\right) \zeta\left(\frac{s}{2}\right)
\zeta\left(\frac{s}{2}-3\right).
\end{equation}
Critical zeros are given by $2^{1-\rho/2}-4^{2-\rho/2}=1$, the off-critical ones 
are localized on the axes $\rho_x=1$ and $\rho_x=7$.

For small odd dimensions only approximate formulas with controlled remainders 
were found \cite{Edery06}.
Explicit formulas for the Epstein zeta-function in terms of the Riemann
and Hurwitz zeta-functions were derived in \cite{Elizalde89b,Kirsten94}.
The Laurent series expansion about the singular point $s=d$ was the
subject of the work \cite{Joyce16}.
The minima and convexity of the Epstein zeta-function were investigated
in \cite{Rankin53,Cassels59,Diananda64,Sarnak06,Lim08}.
General results about the distribution of the Epstein zeros were derived in
\cite{Hejhal87,Potter35,Bombieri87,Steuding05,Nakamura13}.
A fast numerical algorithm for the evaluation of the $d$-dimensional Epstein 
zeta-function in the entire $s$-plane was developed in \cite{Crandall98}.

Critical zeros of the hypercubic Epstein zeta-function are confined to 
the critical line $\Re(s) = \frac{d}{2}$ and it is relatively 
simple to solve numerically one nonlinear equation for 
their imaginary components.
On the other hand, to find blindly the positions of all
off-critical zeros is a hopeless task for dimensions $d\ne 4,6,8$.
It is the aim of this paper to establish a generation mechanism of
off-critical zeros. 
In particular, we apply the critical theory of continuous second-order phase 
transitions in many-body systems, which is developed within the 
condensed-matter and equilibrium statistical physics \cite{Baxter82,Samaj13}, 
to the mathematical problem of generation of Epstein's off-critical zeros 
as bifurcations from specific critical zeros.

To accomplish our aim, we first perform an analytical continuation of 
the Epstein zeta-function to the whole complex $s$-plane, with the spatial 
dimension $d$ being a continuous variable ranging from $0$ to $\infty$.
The consequent numerical results for critical zeros reveal that the latter form 
closed or semi-open curves which enclose disjunctive regions of 
the plane $\left( \rho_x=\frac{d}{2},\rho_y \right)$.
Each curve involves a number of left/right ``edge'' points 
$\rho^*= \left( \frac{d^*}{2},\rho_y^* \right)$, defined by a divergent tangent 
${\rm d}\rho_y/{\rm d}d\vert_{\rho^*}$.
Every edge point gives rise to two conjugate tails of off-critical zeros with 
continuously varying dimension: $d<d^*$ for left and $d>d^*$ for right
edge points.
The curves of critical and off-critical zeros exhibit a singular expansion 
around the critical edge points whose derivation resembles the one around 
a critical point of second-order phase transitions.
The order parameter is identified with the deviation from the critical line:
it is zero along the curve of critical zeros and becomes non-zero along 
the two tails of off-critical zeros.
The singular behavior of the order parameter close to the edge point
is characterized by the mean-field exponent $\frac{1}{2}$.
Various versions of the generation mechanism are discussed.

It turns out that for each $d>9.24555\ldots$ there exists a conjugate 
pair of {\em real} off-critical zeros which tend to 0 and $d$ in the limit 
$d\to\infty$. 
As a by-product of the formalism, we derive an exact result for
$\lim_{d\to 0} \zeta^{(d)}(s)/d$.
An equidistant distribution of critical zeros along the imaginary axis
is obtained for large $d$, spacing between the nearest-neighbour
zeros vanishes as $2\pi/\ln d$ in the limit $d\to\infty$.

The paper is organized as follows.
An analytical continuation of the Epstein zeta-function to the whole complex
$s$-plane and the continuous spatial dimension $d\in (0,\infty)$ 
is constructed in section \ref{Sec2}.
Basic formulas for zeros of the Epstein zeta-function, together with
specific sum rules, are given in section \ref{Sec3}.
The precise location of off-critical zeros of the Epstein zeta-function
in the limit $d\to 0$ and the equidistant distribution of
its critical zeros for large $d$ are derived in section \ref{Sec4}.
Section \ref{Sec5} deals with the numerical evaluation of the curves of 
critical zeros and a singular expansion of these curves around the edge points.
Section \ref{Sec6} describes the generation mechanism of off-critical zeros
from the critical edge points.
Analytical expansion formulas close to the edge points are tested 
numerically.
The concluding section \ref{Sec7} brings a short recapitulation and 
open questions.   

\section{Regularization in $d$ dimensions}\label{Sec2}

\subsection{$\Re(s)>d$}\label{sec2a}
The lattice-sum representation (\ref{zetad}) of $\zeta^{(d)}(s)$ is defined 
for $\Re(s)>d$.
Using the standard Gamma-identity
\begin{equation} \label{Gamma}
\frac{1}{r^s}=\frac{1}{(r^2)^{s/2}} =\frac{1}{\Gamma(s/2)}
\int_{0}^\infty {\rm d} t\, t^{{s/2}-1} {\rm e}^{-r^2 t} ,
\end{equation}
the Epstein zeta-function can be reexpressed as
\begin{eqnarray}
\zeta^{(d)}(s) & = & \frac{1}{2\Gamma(s/2)} \int_{0}^\infty {\rm d} t\, t^{s/2-1} 
\left[ \sum_{n_1,n_2,\ldots, n_d=-\infty}^{\infty} {\rm e}^{-(n_1^2+n_1^2+\ldots+n_d^2)t} - 1 
\right] \nonumber \\ 
& = & \frac{1}{2\Gamma(s/2)} \int_{0}^\infty {\rm d} t\, t^{s/2-1}
\left[ \theta_3^d\left({\rm e}^{-t}\right) - 1 \right] , \label{zetadd}
\end{eqnarray}
where we introduced the Jacobi elliptic function with zero argument  
$\vartheta_3(0,q)\equiv\theta_3(q)=\sum_{n=-\infty}^{\infty} q^{n^2}$ 
(see \cite{Gradshteyn}) and $-1$ subtracts the summand
with $n_1=n_2=\ldots=n_d=0$.
The last integral in (\ref{zetadd}) is known as the Mellin transform of 
the function in the square bracket.

The elliptic theta function $\theta_3\left({\rm e}^{-t}\right)$ exhibits
the following small-$t$ and large-$t$ expansions:
\begin{equation} \label{t3exp}
\theta_3\left({\rm e}^{-t}\right) \mathop{\sim}_{t\to 0}
\sqrt{\frac{\pi}{t}} \left(1 + 2 {\rm e}^{-\pi^2/t}+\cdots \right), \qquad
\theta_3\left({\rm e}^{-t}\right) \mathop{\sim}_{t\to \infty} 
1+2 {\rm e}^{-t} +\cdots . 
\end{equation}
The function under integration in (\ref{zetadd}) is integrable at large $t$ 
and it behaves like $t^{(s-d)/2-1}$ for $t\to 0$, so the real part of 
the power must be greater than $-1$ which yields the mentioned 
restriction $\Re(s)>d$.

To derive another representation of (\ref{zetadd}), we first substitute 
$t=\pi t'$ and then split the integration interval into $0<t'<1$ 
and $1<t'<\infty$, to obtain
\begin{eqnarray} 
\frac{2\Gamma(s/2)}{\pi^{s/2}} \zeta^{(d)}(s) & = &  
\int_0^1 {\rm d} t\, t^{s/2-1} 
\left[ \theta_3^d\left({\rm e}^{-\pi t}\right) - 1 \right] \nonumber \\ 
& & + \int_1^{\infty} {\rm d} t\, t^{s/2-1} 
\left[ \theta_3^d\left({\rm e}^{-\pi t}\right) - 1 \right] . \label{eqq} 
\end{eqnarray}
The Poisson summation formula 
\begin{equation} \label{poisson}
\sum_n {\rm e}^{-(n+\phi)^2 t} =\sqrt{\frac{\pi}{t}}
\sum_n {\rm e}^{2\pi i n\phi} {\rm e}^{-(\pi n)^2/t}
\end{equation}
with $\phi=0$ yields the following relation for the Jacobi theta function 
\begin{equation} \label{PSF}
\theta_3\left({\rm e}^{-\pi t}\right) \equiv \sum_{n=-\infty}^{\infty}
{\rm e}^{-n^2\pi t} = \frac{1}{\sqrt{t}} \sum_{n=-\infty}^{\infty}
{\rm e}^{-n^2\pi/t} = \frac{1}{\sqrt{t}} \theta_3\left({\rm e}^{-\pi/t}\right) . 
\end{equation}
Applying this equality in the last integral of Eq. (\ref{eqq}) and
afterwards using the substitution $t=1/t'$, one finds that 
\begin{equation} \label{eqqq}
\int_1^{\infty} \frac{{\rm d} t}{t} t^{s/2} 
\left[ \theta_3^d\left({\rm e}^{-\pi t}\right) - 1 \right] =
\int_0^1 \frac{{\rm d} t}{t} t^{(d-s)/2} 
\left[ \theta_3^d\left({\rm e}^{-\pi t}\right) - t^{-d/2} \right] .
\end{equation}
As the next step, one adds $-t^{-d/2}+t^{-d/2}$ in the square bracket of 
the first integral on the rhs of (\ref{eqq}) and integrates explicitly 
the remaining term $-1+t^{-d/2}$, which can be done for $\Re(s)>d$.
The final formula reads as
\begin{eqnarray} 
\pi^{-s/2} \Gamma\left(\frac{s}{2}\right) \zeta^{(d)}(s) 
& = & - \frac{1}{s} - \frac{1}{d-s} \nonumber \\ & & + \frac{1}{2} 
\int_0^1 \frac{{\rm d}t}{t} \left( t^{s/2} + t^{(d-s)/2} \right)
\left[ \theta_3^d\left({\rm e}^{-\pi t}\right) - t^{-d/2} \right] .
\phantom{aaa} \label{zetadbest} 
\end{eqnarray}
By using the first relation in (\ref{t3exp}), the difference
$\theta_3^d\left({\rm e}^{-\pi t}\right) - t^{-d/2} \sim 2d
t^{-d/2} {\rm e}^{-\pi/t}$ for $t\to 0$ and the integral on the rhs 
converges for any complex $s$. 
The representation (\ref{zetadbest}) is therefore an analytic continuation 
of (\ref{zetadd}) to the whole complex plane, except for the simple
pole at $s=d$.
The limit $s\to 0$ does not represent any problem as the singularity $-1/s$ 
on the rhs has a counterpart $\Gamma(s/2)\sim_{s\to 0} 2/s$ on the lhs, 
so that $\zeta^{(d)}(0) = - 1/2$ for any $d$.

For $d=1$, the derivation of the formula (\ref{zetadbest}) goes back to 
Riemann's 1859 paper \cite{Riemann1859}
see also an alternative representation (1.12) on page 12 
of the book \cite{Kirsten10}.
For $d>1$, a formula analogous to (\ref{zetadbest}), with the integration
range constrained to $t\in [1,\infty)$, was derived in \cite{Terras80}.

The crucial representation (\ref{zetadbest}) can be used to calculate 
the Epstein zeta-function for any complex $s$ and we could stop here
our analytical analysis.
However, to show the consistency of the formalism with other approaches
(working for $d=1,2,3$) and its relation to neutral Coulomb systems
in thermal equilibrium, in what follows we shall propose other analytic 
continuations of the lattice formula (\ref{zetadd}) to the regions 
$0<\Re(s)<d$ and $\Re(s)<0$ and show that they all lead to 
the same representation (\ref{zetadbest}).

\subsection{$0<\Re(s)<d$} \label{sec2b}
To ensure the convergence of $\zeta^{(d)}(s)$ for $0<\Re(s)<d$, it is useful 
to ''neutralise'' the ``charged'' particle systems by a system of opposite
charge in the way it is often made in Coulomb systems \cite{Samaj12}.

One possibility is to introduce the opposite charges on one half of the sites 
by inserting the factor $(-1)^{n_1}$ or $(-1)^{n_1+n_2}$, etc. into the sum 
(\ref{zetad}) and then express $\zeta^{(d)}(s)$ by using these finite 
expressions \cite{Borwein13}. 
To generate such expressions in a compact form, let us introduce another 
Jacobi theta function $\theta_4(q)=\sum_{n=-\infty}^{\infty} (-1)^nq^{n^2}$.
The function $\theta_4\left( {\rm e}^{-t}\right)$ exhibits the following
small-$t$ and large-$t$ expansions: 
\begin{equation} \label{t4exp}
\theta_4\left({\rm e}^{-t}\right) \mathop{\sim}_{t\to 0}
\sqrt{\frac{\pi}{t}} {\rm e}^{-\pi^2/(4t)}+\cdots , \qquad
\theta_4\left({\rm e}^{-t}\right) \mathop{\sim}_{t\to \infty} 
1-2 {\rm e}^{-t} + \cdots . 
\end{equation}
One of the equalities fulfilled by the Jacobi theta functions reads as 
\cite{Borwein13,Whittacker96}
\begin{equation}\label{ww34}
\theta_3(q)= \frac{1}{2} \left[ \theta_3\left( q^{1/4}\right) 
+ \theta_4\left( q^{1/4}\right) \right] .
\end{equation}
Replacing $\theta^d_3\left( {\rm e}^{-t}\right)$ in (\ref{zetadd})
by using this formula, substituting $t'=t/4$, using the binomial expansion
formula and finally putting ${d\choose 0}\theta^d_3(e^{-t})-1$ on the lhs
of the equation, one arrives at
\begin{eqnarray}
\zeta^{(d)}(s) & = & \frac{1}{2\Gamma(s/2)(2^{d-s}-1)}
\int_{0}^\infty {\rm d} t\ t^{s/2-1} \nonumber \\ & & \times
\left[\sum_{n=1}^{d} {d\choose n} 
\theta_3^{d-n}\left({\rm e}^{-t}\right) \theta_4^n\left({\rm e}^{-t}\right)
-2^d+1\right]. \label{zetabin}
\end{eqnarray}
For large $t$, the expression in the square bracket $\sim -2d{\rm e}^{-t}$
which, when multiplied by $t^{s/2-1}$, is an integrable function.
The small-$t$ asymptotic formula (\ref{t4exp}) for 
$\theta_4\left({\rm e}^{-t}\right)$ ensures that the integral on the rhs 
converges for $\Re(s)>0$. 
This approach was formulated for $d=3$ in \cite{Borwein13}.

Another (more direct) way to regularize the Epstein zeta-function is 
the introduction of a homogeneous neutralising background which cancels
an infinite constant from the hypercubic summation (\ref{zetad})
in the critical strip $0<\Re(s)<d$. 
The regularization procedure is known in Coulomb jellium models; 
for a detailed explanation for the three-dimensional Coulomb
potential $(s=1)$ in the spatial dimension $d=2$, see \cite{Samaj12}. 
To extend the regularization procedure to any $d$ and $s$, let us consider a 
$d$-dimensional sphere of radius $R$ around the reference particle
at the point $(0,\ldots,0)$ and restrict the sum (\ref{zetad}) to 
particles at points with coordinates $n_1^2+n_2^2+\ldots+n_d^2\le R^2$. 
The neutralising background of unit density (equivalent to the
particle density) and opposite ``charge'' sign, localized inside the sphere, 
interacts with the reference particle by the potential
\begin{equation}
-\frac{1}{2} \int_0^R {\rm d}^d r \frac{1}{\vert {\bf r}\vert^s}
= -\frac{1}{2} \int_0^R {\rm d}r\, s_d r^{d-1} \frac{1}{r^s} ,
\end{equation}
where $s_d = 2\pi^{d/2}/\Gamma(d/2)$ is the surface area of the $d$-dimensional
unit sphere. 
The application of the Gamma-identity (\ref{Gamma}) to $1/r^s$ and 
the integration over $r$ results in
\begin{equation}
- \frac{1}{2 \Gamma(s/2)} \int_0^{\infty} {\rm d}t\, 
t^{s/2-1} \left( \frac{\pi}{t} \right)^{d/2}
\frac{1}{\Gamma(d/2)} \left[ \Gamma(d/2) - \Gamma(d/2,R^2 t) \right] ,
\end{equation}
where $\Gamma(d/2,R^2 t)$ is the incomplete Gamma function. 
Inserting the corresponding part of this expression into the square bracket of 
the integrated function in (\ref{zetadd}), $\Gamma(d/2,R^2 t)$ goes to zero as 
the radius $R$ of the $d$-dimensional sphere with the neutralising background 
goes to infinity. 
This leads to the addition of a background term $-(\pi/t)^{d/2}$ in 
the square bracket, i.e.,
\begin{equation} \label{zetad0sd}
\zeta^{(d)}(s) = \frac{1}{2\Gamma(s/2)} \int_{0}^\infty {\rm d} t\, t^{s/2-1}
\left[\theta_3^d\left({\rm e}^{-t}\right) - 1 -
\left(\frac{\pi}{t}\right)^{d/2} \right] 
\end{equation}
with $0<\Re(s)<d$.
The addition of the background term exactly cancels the leading term of 
the expansion of $\theta_3^d\left({\rm e}^{-t}\right)$ at small $t$, 
see the first relation in (\ref{t3exp}), removing in this way the divergence 
of the integral.
The term $-t^{s/2-1}$ is integrable at small $t$ for $\Re(s)>0$.
The term dominant at large $t$ is proportional to $t^{-1+(s-d)/2}$ and it
is integrable for $\Re(s)<d$ as it should be.
The equivalence of the representations (\ref{zetabin}) and (\ref{zetad0sd})
can be proved easily by applying the equality (\ref{ww34}).

Starting from (\ref{zetad0sd}), one can apply the procedure analogous to that 
between equations (\ref{eqq}) and (\ref{eqqq}). 
Integrating explicitly the $-1$ terms in the last step, we arrive at 
the same formula (\ref{zetadbest}) as before, confirming its validity 
also for the region $0<\Re(s)<d$.

\subsection{$\Re(s)\le 0$} \label{sec2c}
We keep the term $n_1=n_2=\ldots=n_d=0$ in the sum in (\ref{zetad}) 
for $\Re(s)<0$, as the summand $r^{-s}$ vanishes automatically for $r=0$.
Omitting $-1$ in the square bracket of (\ref{zetad0sd}), which corresponds to
the cancellation of the term $(0,\ldots,0)$ from the summation, and
maintaining the neutralising background term, one has
\begin{equation} \label{zetads0}
\zeta^{(d)}(s)=\frac{1}{2\Gamma(s/2)} \int_{0}^\infty {\rm d} t\, t^{s/2-1}
\left[\theta_3^d\left({\rm e}^{-t}\right) -\left(\frac{\pi}{t}\right)^{d/2}
\right] , \qquad \Re(s) < 0 .
\end{equation}
While the function under integration is always integrable in the
region of small $t$, its large-$t$ limit $t^{s/2-1}$ is integrable
for $\Re(s)<0$.  
Starting from the representation (\ref{zetads0}) and repeating the steps 
presented in subsection \ref{sec2a}, we recover once more the universal 
relation (\ref{zetadbest}) valid in the whole complex $s$-plane.

\subsection{Functional relation and non-integer values of $d$} \label{sec2d}
Another property of the crucial relation (\ref{zetadbest}) is that its rhs 
is invariant under the transformation $s\to d-s$. 
This symmetry implies the well-known functional relation 
for the self-dual hypercubic lattices \cite{Blanc15}
\begin{equation} \label{funrel}
\pi^{-s/2} \Gamma\left(\frac{s}{2}\right) \zeta^{(d)}(s) =
\pi^{-(d-s)/2} \Gamma\left(\frac{d-s}{2}\right) \zeta^{(d)}(d-s).
\end{equation}
This relation provides an analytical continuation of the lattice sum 
(\ref{zetad}) from the region $\Re(s)>d$ to $\Re(s)<0$ and it can serve as 
another check of the validity of the representation (\ref{zetads0}).

While the original lattice-sum representation of the Epstein zeta-function 
(\ref{zetad}) as well as the representation (\ref{zetabin}) are defined
only for positive integer values of the dimension $d$, $d$ can change 
continuously from $0$ to $\infty$ in the rhs of (\ref{zetadbest}).
In other words, the formula (\ref{zetadbest}) corresponds to an
extension of the definition of the lattice sum (\ref{zetad}) to non-integer
dimensions.   
The possibility to treat (positive) continuous values of $d$ is 
of fundamental importance in this work.
Although the limit $d\to 0^+$ is problematic in the physical sense,
it is well defined algebraically.  

\section{Definition of zeros, sum rules} \label{Sec3}
The ``trivial'' zeros are related to the divergence of the Gamma functions
$\Gamma(s/2)$ at $s=-2n$ $(n=1,2,\ldots)$.
The rhs of Eq. (\ref{zetadbest}) does not exhibit any zeros at these 
trivial points.

The ``critical'' zeros (on the critical line) have $\rho_x=\frac{d}{2}$. 
Inserting $\rho = \frac{d}{2} + {\rm i} \rho_y$ into the rhs of 
(\ref{zetadbest}), the expression becomes real and its nullity determines
the imaginary component $\rho_y$ as follows
\begin{equation} \label{standard}
- \frac{d}{(d/2)^2+\rho_y^2} + \int_0^1 {\rm d}t\, t^{d/4-1}
\cos\left( \frac{\rho_y \ln t}{2} \right) 
\left[ \theta_3^d\left({\rm e}^{-\pi t}\right) - t^{-d/2} \right] = 0 .   
\end{equation}
This equation is symmetric with respect to the complex conjugation
$\rho_y\to -\rho_y$.
In $d=1$, the critical zeros with $\rho_x=1/2$ are those suggested
by Riemann to be the only nontrivial ones in the complex plane.

The ``off-critical'' zeros (off the critical line) are those with 
$\rho_x\ne \frac{d}{2}$. 
Let us denote the deviation of $\rho_x$ from its critical value $\frac{d}{2}$ as
\begin{equation} \label{Deltarho}
\Delta\rho_x(d) = \rho_x(d) - \frac{d}{2} . 
\end{equation}
In this case, the rhs of (\ref{zetadbest}) becomes complex and the off-critical 
zeros are given by the pair of coupled equations 
\begin{eqnarray}
- \left[ \frac{\frac{d}{2}+\Delta\rho_x}{\left(\frac{d}{2}+\Delta\rho_x\right)^2
+\rho_y^2} + \frac{\frac{d}{2}-\Delta\rho_x}{
\left(\frac{d}{2}-\Delta\rho_x\right)^2+\rho_y^2} \right]  
\phantom{aaaaaaaaaaaaaaaaaaaa} \nonumber \\
+ \int_0^1 {\rm d}t\, t^{d/4-1} \cos\left( \frac{\rho_y \ln t}{2} \right) 
\cosh\left( \frac{\Delta\rho_x \ln t}{2} \right) 
\left[ \theta_3^d\left({\rm e}^{-\pi t}\right) - t^{-d/2} \right] = 0,
\label{nonstandard1} \\   
\rho_y \left[ \frac{1}{\left(\frac{d}{2}+\Delta\rho_x\right)^2+\rho_y^2} 
- \frac{1}{\left(\frac{d}{2}-\Delta\rho_x\right)^2+\rho_y^2} \right]  
\phantom{aaaaaaaaaaaaaaaaaaa} \nonumber \\
+ \int_0^1 {\rm d}t\, t^{d/4-1} \sin\left( \frac{\rho_y \ln t}{2} \right) 
\sinh\left( \frac{\Delta\rho_x \ln t}{2} \right) 
\left[ \theta_3^d\left({\rm e}^{-\pi t}\right) - t^{-d/2} \right] = 0 . 
\label{nonstandard2}
\end{eqnarray}
These equations are symmetric with respect to the sign reversals
$\Delta\rho_x\to -\Delta\rho_x$ and $\rho_y\to -\rho_y$.
This means that when $(\Delta\rho_x,\rho_y)$ with $\Delta\rho_x\ne 0$ is 
the zero solution of equations (\ref{nonstandard1}) and (\ref{nonstandard2}),
also $(\Delta\rho_x,-\rho_y)$, $(-\Delta\rho_x,\rho_y)$ and
$(-\Delta\rho_x,-\rho_y)$ belong to the set of zero points.  

The critical and off-critical zeros satisfy certain constraints (sum rules)
which follow from the universal representation (\ref{zetadbest}).
Let us rewrite that representation as
\begin{eqnarray} 
h^{(d)}(s) & \equiv & \frac{s(s-d)}{d} \pi^{-s/2} \Gamma(s/2) \zeta^{(d)}(s)
\nonumber \\
& = & 1 + \frac{s(s-d)}{2d} \int_0^1 \frac{{\rm d}t}{t} 
\left( t^{s/2} + t^{(d-s)/2} \right) \left[ \theta_3^d\left({\rm e}^{-\pi t}\right)
- t^{-d/2} \right] . \label{first} 
\end{eqnarray}
$h^{(d)}(s)$ is an entire function of $s$ such that $h^{(d)}(0)=1$.
It vanishes at the nontrivial (critical and off-critical) zeros
$\{ \rho_n \}_{n=1}^{\infty}$ of $\zeta^{(d)}(s)$ (to simplify notation,
we omit the upper index $(d)$ in $\rho_n$). 
According to the Weierstrass factorization theorem, 
$h^{(d)}(s)$ can be factored over its nontrivial zeros
$\{ \rho_n\}_{n=1}^{\infty}$.
Let $p(d)$ be the smallest non-negative integer such that the series
\begin{equation} \label{pd}
\sum_{n=1}^{\infty} \frac{1}{\vert \rho_n\vert^{p(d)+1}}
\end{equation}
converges.
Then $h^{(d)}(s)$ is expressible in terms of the Hadamard's product
\begin{equation} \label{Hadamard1} 
h^{(d)}(s) = {\rm e}^{-\sum_{n=1}^{\infty}\sum_{k=1}^{p(d)}\frac{s^k}{k \rho_n^k}} 
\prod_{n=1}^{\infty} \left( 1 - \frac{s}{\rho_n} \right)
{\rm e}^{\sum_{k=1}^{p(d)} \frac{s^k}{k \rho_n^k}} ,      
\end{equation}
where the exponential factor inside the product over $n$ ensures the product
convergence. 
The corresponding sum rules for the inverse powers of zeros follow
from the generating formula \cite{McPhedran18}
\begin{equation} \label{srd}
\sigma_k \equiv \sum_{n=1}^{\infty} \frac{1}{\rho^k} 
= - \frac{1}{(k-1)!} \frac{{\rm d}^k}{{\rm d}s^k}
\ln h^{(d)}(s) \Bigg\vert_{s=0} .
\end{equation}
This result was derived by comparing two expansions, one based on the Taylor
series for $\ln [h^{(d)}(s)/h^{(d)}(0)]$ and the other based on the product
representation (\ref{Hadamard1}), see Eqs. (4)-(6) of \cite{McPhedran18}. 
As was mentioned in \cite{McPhedran18}, the same result can be obtained
formally by taking the product representation of $h^{(d)}(s)$ without
any convergence factors, i.e.
\begin{equation}
h^{(d)}(s) = \prod_{n=1}^{\infty} \left( 1 - \frac{s}{\rho_n} \right) ,
\end{equation}
and performing the Taylor series of its logarithm around $s=0$, 
\begin{equation}
\ln h^{(d)}(s) = \sum_{n=1}^{\infty} \ln\left( 1 - \frac{s}{\rho_n} \right)
= \sum_{k=1}^{\infty} \frac{s^k}{k} \sum_{n=1}^{\infty} \frac{1}{\rho_n^k} ,   
\end{equation}
which is consistent with the generating formula (\ref{srd}).
In other words, the formula for the zeros moments (\ref{srd}) is
independent of the particular form of the product regularization of $h^{(d)}(s)$.
Inserting the representation (\ref{first}) of $h^{(d)}(s)$ into
the generating formula (\ref{srd}), we obtain the following sum rules
for the inverse powers of zeros:
\begin{eqnarray}
\sum_{\rho} \frac{1}{\rho} & = & \frac{1}{2} 
\int_0^1 \frac{{\rm d}t}{t} (1+t^{d/2}) 
\left[ \theta_3^d\left({\rm e}^{-\pi t}\right) - t^{-d/2} \right] , 
\label{sr1} \\
\sum_{\rho} \frac{1}{\rho^2} & = & \left( \sum_{\rho} \frac{1}{\rho} \right)^2
- \frac{1}{2d} \int_0^1 \frac{{\rm d}t}{t} \left[ 2 (1+t^{d/2}) - d (1-t^{d/2})
\ln t \right] \nonumber \\ & & \times
\left[ \theta_3^d\left({\rm e}^{-\pi t}\right) - t^{-d/2} \right] , \label{sr2} \\
\sum_{\rho} \frac{1}{\rho^3} & = & - \frac{1}{2} 
\left( \sum_{\rho} \frac{1}{\rho} \right)^3 + 
\frac{3}{2} \left( \sum_{\rho} \frac{1}{\rho} \right)
\left( \sum_{\rho} \frac{1}{\rho^2} \right) \nonumber \\ & &
- \frac{3}{16 d} \int_0^1 \frac{{\rm d}t}{t} \ln t
\left[ 4 (1-t^{d/2}) - d (1+t^{d/2})\ln t \right]
\nonumber \\ & & \times 
\left[ \theta_3^d\left({\rm e}^{-\pi t}\right) - t^{-d/2} \right] , \label{sr3} 
\end{eqnarray}
etc. 

In the $d=1$ case of the Riemann zeta-function, $p(1)=1$ and $h^{(1)}(s)$
can be factorized as \cite{Hadamard93}
\begin{equation} \label{Hadamard} 
h^{(1)}(s) = {\rm e}^{[\ln(4\pi)-2-\gamma_0]s/2} 
\prod_{n=1}^{\infty} \left( 1 - \frac{s}{\rho_n} \right) {\rm e}^{s/\rho_n} ,
\end{equation}
where $\gamma_0=0.57721\ldots$ is the Euler-Mascheroni constant.
The corresponding sum rules for the inverse powers of zeros,
derived in \cite{Lehmer88,Keiper92,Li97} (see also p. 168 of \cite{Finch03}), 
read as
\begin{eqnarray}
\sum_{\rho} \frac{1}{\rho} & = & \frac{1}{2} \left[ 2 + \gamma_0 - 
\ln(4\pi) \right] , \nonumber \\
\sum_{\rho} \frac{1}{\rho^2} & = & 1 + \gamma_0^2 - \frac{1}{8}\pi^2
+ 2 \gamma_1 , \nonumber \\
\sum_{\rho} \frac{1}{\rho^3} & = & 1 + \gamma_0^3 + 3 \gamma_0 \gamma_1
+ \frac{3}{2} \gamma_2 - \frac{7}{8} \zeta(3) ,  \label{St} 
\end{eqnarray}
etc.
Here, $\{ \gamma_n\}_{n=1}^{\infty}$ are the Stieltjes constants 
defined by the Laurent expansion of $\zeta(s)$ around the singular point
$s=1$ \cite{Coffey06}:
\begin{equation} \label{Stj} 
\zeta(s) = \frac{1}{s-1} + \sum_{n=0}^{\infty} \frac{(-1)^n}{n!} \gamma_n
(s-1)^n .
\end{equation}
The sum rules (\ref{St}) follow directly from (\ref{sr1})-(\ref{sr3}) 
by taking $d=1$ and expressing integrals over $t$ in terms of special
constants.
Note that the sum $\sum_{\rho} 1/\rho$ is not absolutely convergent.
However, every zero $\rho$ has its complex conjugate $\rho^*$ and
the above sum has to be understood in the following sense:
$\sum_{\rho;\rho_y>0} \left( 1/\rho + 1/\rho^* \right)$;
the summands are proportional to $1/\rho_y^2$ for large $\rho_y$
and therefore this sum is absolutely convergent.
We apply this convention in what follows, whenever 
the sum $\sum_{\rho} 1/\rho$ appears. 

As concerns dimensions $d>1$, to our knowledge there is no general analysis 
in the mathematical literature about the dependence of the smallest
non-negative integer $p$, which ensures the convergence of the series
(\ref{pd}), on the spatial dimension $d$.
Assuming that this integer is finite, one can calculate the $\sigma$-moments
by using Eqs. (\ref{sr1})-(\ref{sr3}).
Let us make a comparison of the relations (\ref{sr1})-(\ref{sr3}) with
the numerical evaluation of these sum rules by taking a large number of
zeros.
It is relatively simple to generate zeros of the Epstein zeta-function for
dimensions $d=2,4,6,8$ by using the corresponding exact formulas
(\ref{zeta2})-(\ref{zeta8}).
For $d=2$, see formula (\ref{zeta2}), in the numerical calculation of 
the inverse powers of zeros with $\Re(\rho)=1$ we took into account 
all zeros of $\zeta(s/2)$ and the first 200 pairs of complex conjugate zeros
of $\beta(s/2)$. 
We got the numerical result $\sigma_1^{\rm num} = 0.0487637\ldots$ 
compared with the exact result (\ref{sr1})
$\sigma_1^{\rm exact} = 0.0504398\ldots$,
$\sigma_2^{\rm num} = - 0.048525\ldots$ compared with (\ref{sr2})
$\sigma_2^{\rm exact} = - 0.050201\ldots$ and
$\sigma_3^{\rm num} = - 0.000355375\ldots$ compared with (\ref{sr3})
$\sigma_3^{\rm exact} = - 0.000355379\ldots$.
For $d=4$, see formula (\ref{zeta4}), we took into account all zeros implied by
the factors $(1-2^{2-s})$ and $\zeta(s/2)$, and the first 10000 pairs of complex
conjugate zeros of $\zeta(s/2-1)$.
We got $\sigma_1^{\rm num} = 0.123502\ldots$ compared with 
$\sigma_1^{\rm exact} = 0.123704\ldots$,
$\sigma_2^{\rm num} = - 0.059365\ldots$ compared with
$\sigma_2^{\rm exact} = - 0.059433\ldots$ and
$\sigma_3^{\rm num} = - 0.001712059732\ldots$ compared with 
$\sigma_3^{\rm exact} = - 0.001712059733\ldots$.
The agreement between the numerical and exact results is very good and
improves itself with increasing the inverse power of zeros, as it should be.

\section{Special limits of dimension $d$} \label{Sec4}
In this section, we consider two special limits of the spatial dimension: 
$d\to 0$ and $d\to\infty$.
Although these dimensional limits look artificial from a practical 
(physical) point of view, they are well defined algebraically and can be used 
to check the accuracy of numerical results.
 
\subsection{$d\to 0^+$} \label{de0}
We study the small-$d$ limit within the representation (\ref{zetadbest}),
rewritten by using the relation between the Jacobi theta functions
(\ref{PSF}) as follows:
\begin{eqnarray} 
2 \Gamma\left(\frac{s}{2}\right) \zeta^{(d)}(s) 
& = & 2 \pi^{s/2} \left( \frac{1}{s-d} - \frac{1}{s} \right)
\nonumber \\ & & + \pi^{s/2} \int_0^1 \frac{{\rm d}t}{t} 
\left( t^{(s-d)/2} + t^{-s/2} \right) 
\left[ \theta_3^d\left({\rm e}^{-\pi/t}\right) - 1 \right] . \label{zetadbest1} 
\end{eqnarray}
Let us apply the limit $d\to 0^+$ to both sides of this equation.
The integrand converges uniformly on the interval $t\in [0,1]$ and so 
the limit $d\to 0^+$ can be moved inside the integral.
Moreover, the function $\theta_3({\rm e}^{-\pi/t})$ is finite for any
$t\in[0,1]$ and we can set 
$\theta_3^d\left({\rm e}^{-\pi/t}\right)-1\sim d 
\ln\theta_3\left({\rm e}^{-\pi/t}\right)$ in the limit $d\to 0^+$. 
Thus,
\begin{equation} \label{zetadbest2} 
2 \Gamma\left(\frac{s}{2}\right) \lim_{d\to 0^+} \frac{\zeta^{(d)}(s)}{d} 
= \frac{2 \pi^{s/2}}{s^2} + \pi^{s/2} \int_0^1 \frac{{\rm d}t}{t} 
\left( t^{s/2} + t^{-s/2} \right) 
\ln\left[ \theta_3\left({\rm e}^{-\pi/t}\right) \right] . 
\end{equation}
The integral on the rhs of (\ref{zetadbest2}) is the sum of two integrals.
The first integral can be expressed by using the relation (\ref{PSF}) 
as follows  
\begin{equation} \label{i1}
\pi^{s/2} \int_0^1 \frac{{\rm d}t}{t} t^{s/2} 
\ln\left[ \theta_3\left({\rm e}^{-\pi/t}\right) \right] 
= - \frac{2 \pi^{s/2}}{s^2} + \int_0^{\pi} \frac{{\rm d}u}{u} u^{s/2} 
\ln\left[ \theta_3\left({\rm e}^{-u}\right) \right] ,
\end{equation}
where we used the substitution $u=\pi t$.
The second integral is expressible as
\begin{equation} \label{i2}
\pi^{s/2} \int_0^1 \frac{{\rm d}t}{t} t^{-s/2} 
\ln\left[ \theta_3\left({\rm e}^{-\pi/t}\right) \right] 
= \int_{\pi}^{\infty} \frac{{\rm d}v}{v} v^{s/2} 
\ln\left[ \theta_3\left({\rm e}^{-v}\right) \right] ,
\end{equation}
with the substitution $v=\pi/t$.
Finally, taking the integral on the rhs of (\ref{zetadbest2}) as
the sum of (\ref{i1}) and (\ref{i2}), one ends up with the result 
\begin{equation} \label{d0}
f(s) \equiv \lim_{d\to 0^+} \frac{\zeta^{(d)}(s)}{d}
= \frac{1}{2\Gamma(s/2)} \int_0^{\infty} {\rm d}t\, t^{s/2-1}
\ln \left[\theta_3\left({\rm e}^{-t}\right)\right] .   
\end{equation}
Considering the asymptotic behavior of $\theta_3\left({\rm e}^{-t}\right)$
in the limits of small and large $t$, this representation is valid 
for $\Re(s)>0$.
 
To express the integral in (\ref{d0}) in terms of standard functions,
the known product representation of the Jacobi theta function \cite{Gradshteyn}
\begin{equation} \label{theta1}
\theta_3(q) = \prod_{j=1}^{\infty} (1-q^{2j}) (1+q^{2j-1})^2
\end{equation}
is not sufficient for our purposes and must be symmetrized.
Let us consider the function
\begin{equation}
g(q) = \prod_{j=1}^{\infty} (1+q^j) (1-q^{2j-1}) = 1 + a_1 q + a_2 q^2 + \cdots
\end{equation}
which is evidently analytic in $q$ around $q=0$.
Since it holds
\begin{equation}
g(q^2) = \prod_{j=1}^{\infty} (1+q^{2j}) (1-q^{2j-1}) (1+q^{2j-1})
= \prod_{j=1}^{\infty} (1+q^j) (1-q^{2j-1}) = g(q) ,
\end{equation}
we have $g(q)=g(q^2)=g(q^4)=\ldots=1$.
The division of the representation (\ref{theta1}) of $\theta_3(q)$ by $g(q^2)$ 
results in
\begin{equation} \label{theta2}
\theta_3(q) = \prod_{j=1}^{\infty} 
\frac{(1-q^{2j})(1+q^{2j-1})}{(1+q^{2j})(1-q^{2j-1})} .
\end{equation}
Consequently,
\begin{equation}
\ln\theta_3(q) = \sum_{j=1}^{\infty} (-1)^j 
\ln \left( \frac{1-q^j}{1+q^j} \right)
= -2 \sum_{j=1}^{\infty} (-1)^j \sum_{k=0}^{\infty} \frac{1}{2k+1} q^{j(2k+1)} .
\end{equation}
Inserting this expansion into (\ref{d0}) yields
\begin{equation}
f(s) = \sum_{j=1}^{\infty} (-1)^{j+1} \frac{1}{j^{s/2}}
\sum_{k=1}^{\infty} \frac{1}{(2k+1)^{s/2+1}} .
\end{equation}
For $\Re(s)>0$, it is easy to derive the following relations
\begin{eqnarray}
\sum_{j=1}^{\infty} (-1)^{j+1} \frac{1}{j^{s/2}} & = & (1-2^{-s/2+1})
\zeta\left( \frac{s}{2} \right) , \nonumber \\
\sum_{k=1}^{\infty} \frac{1}{(2k+1)^{s/2+1}} & = & (1-2^{-s/2-1})
\zeta\left( \frac{s}{2}+1 \right) . 
\end{eqnarray}
Consequently,
\begin{equation} \label{fr}
\lim_{d\to 0^+} \frac{\zeta^{(d)}(s)}{d} = 
(1-2^{-s/2+1}) (1-2^{-s/2-1})
\zeta\left( \frac{s}{2} \right) \zeta\left( \frac{s}{2}+1 \right) .
\end{equation}
This exact result complements the known formulas (\ref{zeta2})-(\ref{zeta8})
for $d=2,4,6,8$.

The zeros of the $d\to 0^+$ Epstein zeta-function are of two kinds.
One of the first two brackets on the rhs of (\ref{fr}) vanishes for
\begin{equation} \label{sk}
s_k^{\pm} = 2 \left( \pm 1 + \frac{2\pi{\rm i}}{\ln 2} k \right) , \qquad
k = \pm 1, \pm 2, \ldots .  
\end{equation}
Note that $s_k^-=-2$ with $k=0$ is a trivial zero and $s_k^+=2$ 
with $k=0$ is not a zero because the limit
$\lim_{s\to2} (1-2^{-s/2+1}) \zeta(s/2) = \ln 2$ is finite.
The second set of zeros has the origin in the product of
the Riemann zeta-functions $\zeta(s/2)\zeta(s/2+1)$.
Let us denote by $\zeta_n = \frac{1}{2} + {\rm i}\rho_y(n)$ 
$(n=1,2,\ldots)$ the successive series of (critical) zeros of $\zeta(s)$ 
in the upper half-plane with $\rho_y(n)>0$ and by 
$\zeta_n = \frac{1}{2} + {\rm i}\rho_y(n)$ $(n=-1,-2,\ldots)$ 
the complex conjugate ones in the lower half-plane with 
$\rho_y(n)=-\rho_y(-n)<0$. 
The zeros of the product $\zeta(s/2)\zeta(s/2+1)$ are given by
\begin{equation}
s_n^{(1)} = 2\zeta_n , \qquad s_n^{(2)} = 2(\zeta_n-1) .
\end{equation}
According to the Riemann hypothesis, they are constrained to the axes 
$\Re(s)=\pm 1$.

Let us check whether the above zeros fulfill the sum rules 
(\ref{sr1})-(\ref{sr3}) taken in the limit $d\to 0^+$.
The $d\to 0^+$ limit of the sum rule (\ref{sr1}) yields
$\sum_{\rho} 1/\rho = 0$.
This result is confirmed by the explicit calculation
\begin{equation}
\sum_{\rho} \frac{1}{\rho} = \frac{1}{2} \sum_n \left(
\frac{1}{\zeta_n} + \frac{1}{\zeta_n-1} \right) + \frac{1}{2} \sum_{k\ne 0}
\left( \frac{1}{1-\frac{2\pi{\rm i}}{\ln 2}k} -
\frac{1}{1+\frac{2\pi{\rm i}}{\ln 2}k} \right) = 0 ,
\end{equation}
where we have used that $1-\zeta_n=\zeta_{-n}$ (provided that the Riemann
hypothesis holds) and summands are coupled as complex-conjugate pairs.

Now let us take the limit $d\to 0^+$ of the second sum rule (\ref{sr2}).
The convergence of the integrand is uniform and therefore
the limit-integral interchange can be done.
With respect to the relation (\ref{PSF}), the difference
$\theta_3^d({\rm e}^{-\pi t}) - t^{-d/2}$ can be expressed as
$t^{-d/2} \left[ \theta_3^d({\rm e}^{-\pi/t}) - 1 \right]$.
Since the positive function $\theta_3({\rm e}^{-\pi/t})$ is finite for
$t\in [0,1]$, we can expand
$\theta_3^d({\rm e}^{-\pi/t}) - 1 \sim d
\ln\left[ \theta_3({\rm e}^{-\pi/t}) \right]$.
Switching back to $\theta_3({\rm e}^{-\pi t})$ by using (\ref{PSF}),
one obtains
\begin{equation} \label{aa}
\sum_{\rho} \frac{1}{\rho^2} = -2 \int_0^1 \frac{{\rm d}t}{t}
\ln \left[ \sqrt{t} \theta_3\left({\rm e}^{-\pi t}\right) \right] .
\end{equation}
With regard to the definition (\ref{d0}) of the function $f(s)$,
one can write
\begin{equation} \label{twoint}
2 \pi^{-s/2} \Gamma\left(\frac{s}{2}\right) f(s) =
\int_0^1 \frac{{\rm d}t}{t} t^{s/2}
\ln \left[ \theta_3\left({\rm e}^{-\pi t}\right) \right] 
+ \int_1^{\infty} \frac{{\rm d}t}{t} t^{s/2}
\ln \left[ \theta_3\left({\rm e}^{-\pi t}\right) \right] . 
\end{equation}
As concerns the second integral on the rhs, applying the relation (\ref{PSF})
and subsequently using the substitution $t=1/x$ results in
\begin{equation}
\int_0^1 \frac{{\rm d}x}{x} x^{-s/2}
\ln \left[ \sqrt{x} \theta_3\left({\rm e}^{-\pi x}\right) \right] . 
\end{equation}
Adding the ``missing'' prefactor $\sqrt{t}$ to $\theta_3$ in the first
integral on the rhs of (\ref{twoint}) and evaluating the integral
\begin{equation}
- \frac{1}{2} \int_0^1 \frac{{\rm d}t}{t} t^{s/2} \ln t
= \frac{2}{s^2} ,  
\end{equation}  
we end up with
\begin{equation}
2 \pi^{-s/2} \Gamma\left(\frac{s}{2}\right) f(s) = \frac{2}{s^2}
+  \int_0^1 \frac{{\rm d}t}{t} \left( t^{s/2} + t^{-s/2} \right) 
\ln \left[ \sqrt{t} \theta_3\left({\rm e}^{-\pi t}\right) \right] .
\end{equation}
Applying to both sides the limit $s\to 0$ and noting the
uniform convergence of the integrand, with regard to (\ref{aa})
it holds that
\begin{eqnarray}
\sum_{\rho} \frac{1}{\rho^2} & = & -2 \lim_{s\to 0}
\left[ \pi^{-s/2} \Gamma\left( \frac{s}{2} \right) f(s)
- \frac{1}{s^2} \right] \nonumber \\
& = & \frac{\gamma_0^2}{2} - \frac{\pi^2}{16} + (\ln 2)^2 + \gamma_1 .
\end{eqnarray} 
Here, the last equality results from the expansion of 
$\pi^{-s/2}\Gamma(s/2)f(s)$, with $f(s)$ defined by (\ref{d0}) and (\ref{fr}), 
around $s=0$ by using \textit{Mathematica}.
On the other side, summing 
\begin{equation}
\frac{1}{4} \sum_n\left( \frac{1}{\zeta_n^2} + \frac{1}{(\zeta_n-1)^2}\right) 
= \frac{1}{2} \sum_n \frac{1}{\zeta_n^2} = \frac{\gamma_0^2}{2} + \frac{1}{2} 
- \frac{\pi^2}{16} + \gamma_1 ,
\end{equation}
see section 2.21 of \cite{Finch03}, and
\begin{equation}
\frac{1}{4} \sum_{k\ne 0} \left[ 
\frac{1}{\left(1-\frac{2\pi{\rm i}}{\ln 2}k\right)^2} 
+\frac{1}{\left(1+\frac{2\pi{\rm i}}{\ln 2}k\right)^2} \right] =  
-\frac{1}{2} + (\ln 2)^2 ,
\end{equation}
the same result holds by substituting directly the spectrum of zeros
into the sum $\sum_{\rho} 1/\rho^2$.

Finally, from (\ref{sr3}) one gets that $\sum_{\rho} 1/\rho^3 = 0$
which is trivially reproduced by the identified zeros.   

The critical zeros with $\rho_x=0$ are absent at zero dimension.
We conclude that all zeros are off-critical in the limit $d\to 0^+$
and they are constrained to the axes $\Re(s)=\pm 1,\pm 2$. 

\subsection{$d\to\infty$} \label{deinfinity}
By using the relation (\ref{PSF}), the difference
$\theta_3^d\left({\rm e}^{-\pi t}\right) - t^{-d/2}$ can be expanded as
\begin{eqnarray} 
\theta_3^d\left({\rm e}^{-\pi t}\right) - t^{-d/2}
& = & \frac{2d}{t^{d/2}} \Bigg[ {\rm e}^{-\pi/t} + (d-1) {\rm e}^{-2\pi/t}
\nonumber \\ & & + \frac{2}{3} (d-1) (d-2) {\rm e}^{-3\pi/t} + \cdots \Bigg] .
\label{dinfinity}
\end{eqnarray} 
Let us look for the critical zeros in the limit $d\to\infty$. 
Inserting this expansion into (\ref{standard}) and taking the limit 
$d\to\infty$, one gets the condition
\begin{eqnarray} 
\int_0^1 \frac{{\rm d}t}{t} t^{-d/4} \cos\left( \frac{\rho_y \ln t}{2} \right) 
\left[ {\rm e}^{-\pi/t} + (d-1) {\rm e}^{-2\pi/t} \right. & & \nonumber \\
\left.
+ \frac{2}{3} (d-1) (d-2) {\rm e}^{-3\pi/t} + \cdots \right] & = & 0 .   
\label{infini}
\end{eqnarray}
The leading term (in the limit $d\to\infty$) inside the square bracket 
is the first one with ${\rm e}^{-\pi/t}$.
Like for instance, the second term with ${\rm e}^{-2\pi/t}$ gives, after 
the substitution $t=2t'$, a vanishing contribution of order $2^{-d/4} d$ in 
comparison with the first one.
The consequent equation
\begin{equation} \label{rov}
\int_0^1 {\rm d}t\, {\rm e}^{S(t)}
\cos\left( \frac{\rho_y \ln t}{2} \right) = 0 , \qquad
S(t) = -\left( 1 + \frac{d}{4}\right) \ln t - \frac{\pi}{t}
\end{equation} 
can be treated within the Laplace integral method.
The maximum point of the ``action'' $S(t)$, given by the condition 
$\partial S(t)/\partial t\vert_{t_0} = 0$, is $t_0 = 4\pi/d$ in the 
large-$d$ limit.
This point lies in the interval $[0,1]$ as is needed.
The action can be expanded around $t_0$ into a Taylor series as follows 
\begin{equation}
S(t) = S(t_0) - \frac{\pi}{2} \left( \frac{d}{4\pi} \right)^3 (t-t_0)^2
+ \cdots . 
\end{equation}
Considering this expansion in (\ref{rov}), in the limit $d\to\infty$
the exponential ${\rm e}^{S(t)}$ becomes proportional to the Dirac delta 
function $\delta(t-t_0)$, so that the possible values of $\rho_y$ are 
determined by
\begin{equation} \label{infinity}
\cos\left[ \frac{\rho_y}{2}\ln\left( \frac{4\pi}{d}\right) \right] = 0 ,
\qquad \rho_y = \frac{(2n+1)\pi}{\ln[d/(4\pi)]} \quad (n=0,\pm 1,\ldots).
\end{equation}
The critical zeros are therefore distributed equidistantly along 
the imaginary axis in the large-$d$ limit.
Of course, there are next (subleading) terms in the expansion of $\rho_y$
which vanish in comparison with the leading term (\ref{infinity})
for sufficiently large $d$. 
The threshold dimension $d$ beyond which (\ref{infinity}) describes adequately
the distribution of first (say $3, 4,\ldots$) zeros can be found
by comparing with numerical data.  
Since the distance between the nearest-neighbour zeros along the imaginary axis
is proportional to $1/\ln d$ for large values of $d$, the critical zeros
collapse into one point on the real axis (with coordinate $\rho_x=d/2$)
in the limit $d\to\infty$.

\section{Singular expansion around critical edge zeros} \label{Sec5}
\begin{figure}[th]
\centering
\includegraphics[clip,width=0.75\textwidth]{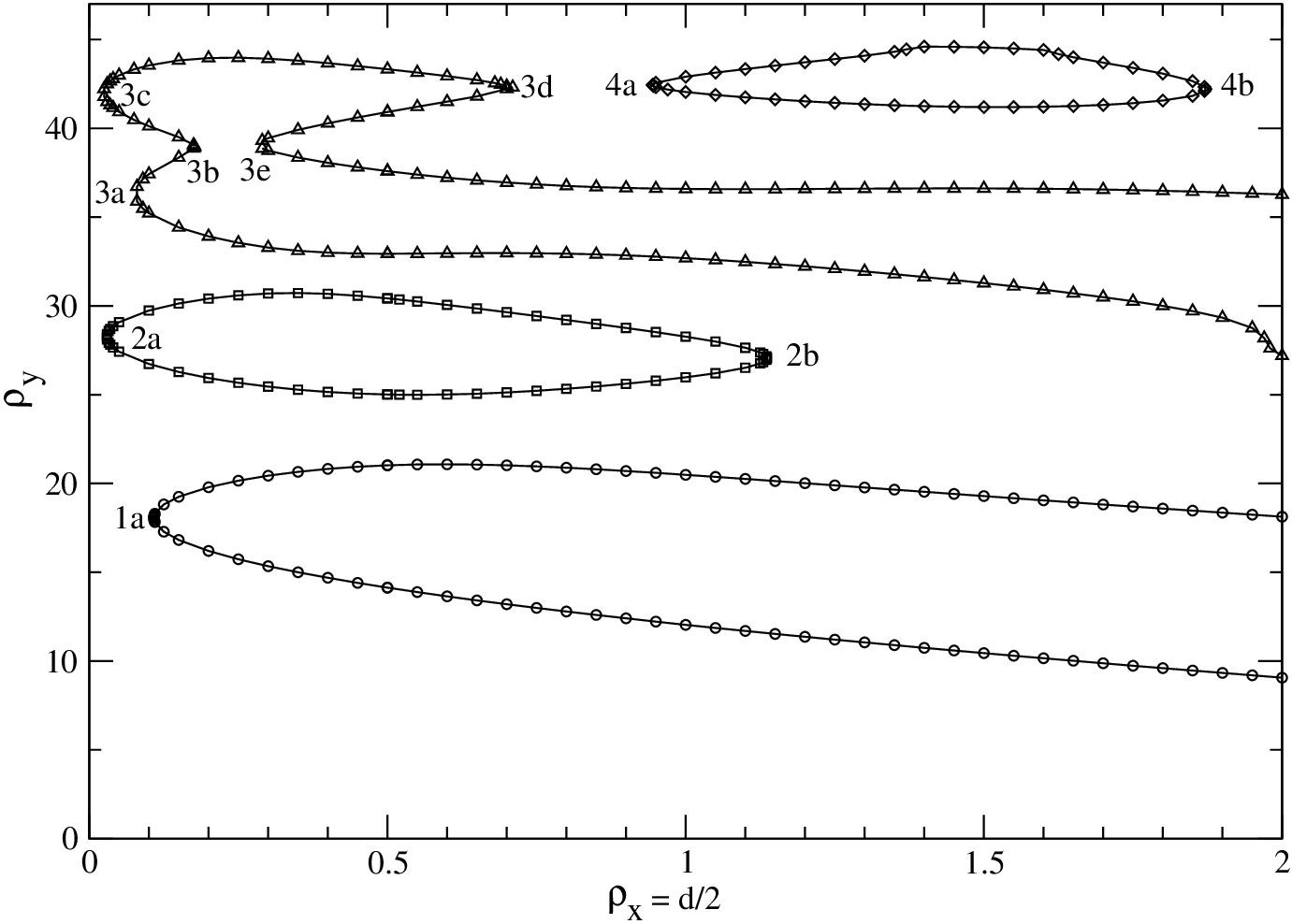}
\caption{Numerical results for the critical Epstein zeros at small dimensions 
$d\le 4$. 
The critical zeros are connected by solid curves, closed or semi-open, 
which enclose disjunctive regions of the complex plane.
The left and right edge points of the curves are indicated as 
1a, 2a, 3a, 3c, 3e, 4a and 2b, 3b, 3d, 4b, respectively.}
\label{alldpol}
\end{figure}

\begin{figure}[th]
\centering
\includegraphics[clip,width=0.75\textwidth]{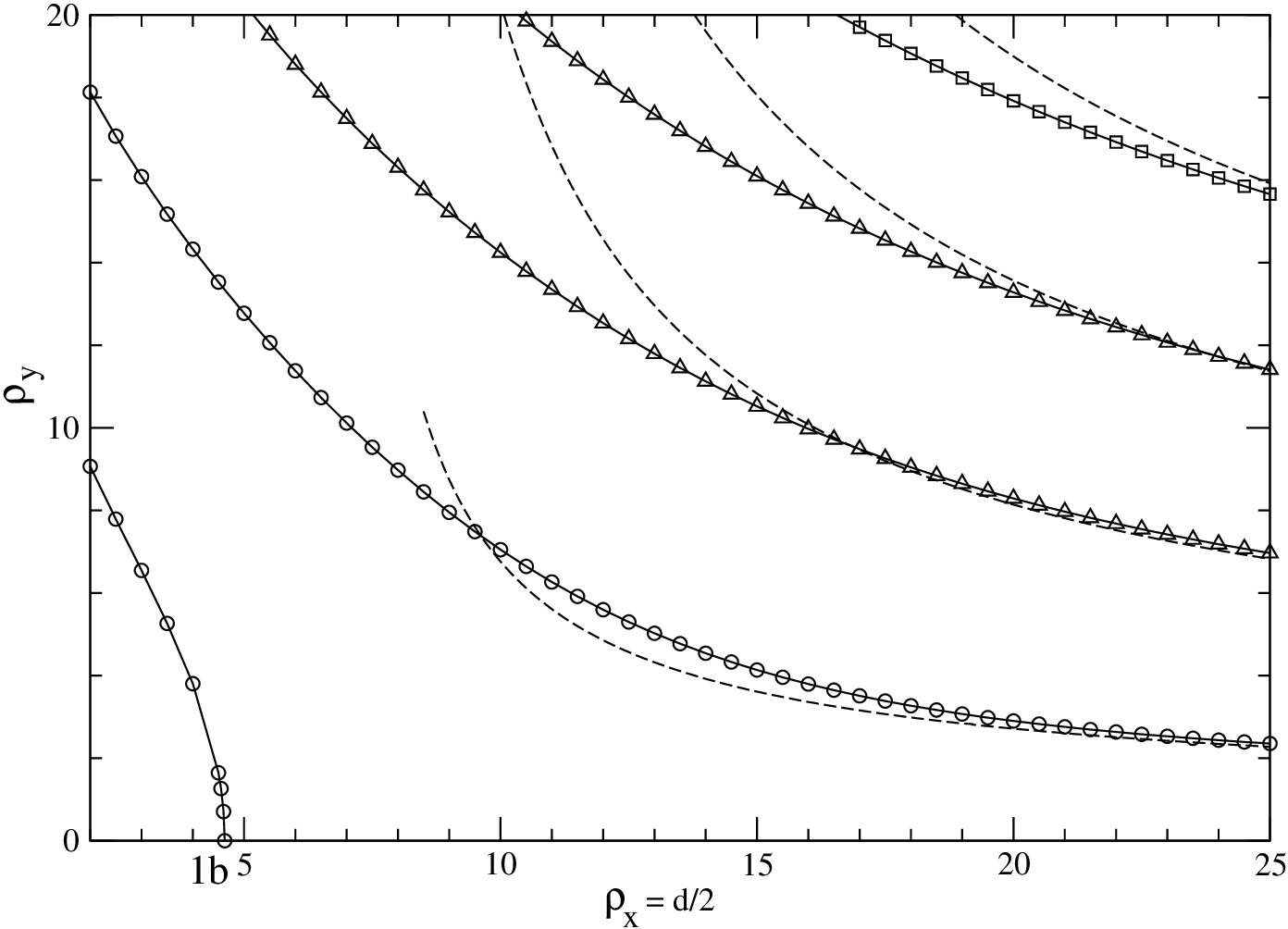}
\caption{Numerical results for critical Epstein zeros at large dimensions
(up to $d=50$).
The lowest curve ending up at the edge point 1b corresponds to 
the prolongation of the lowest curve in Fig. \ref{alldpol}.
The solid lines of the second (open circles) and third (open triangles) zeros 
(as they appear in the region of small $d$) go up to $d\to\infty$.
The prediction of zeros following from the asymptotic relation (\ref{infinity})
is represented by the dashed lines.}
\label{larged3}
\end{figure}

In Figs. \ref{alldpol} and \ref{larged3}, we present by open symbols 
interconnected via solid lines the numerical results for the critical zeros of 
the Epstein zeta-function $\zeta^{(d)}(s)$, with the imaginary components 
$\rho_y$ smaller than 45 to keep high numerical accuracy (6-25 decimal digits)
of their determination.
For a fixed dimension $d$, or $\rho_x=d/2$, the numerical evaluation of 
one zero using \textit{Mathematica} takes approximately 5 seconds of CPU time 
on a standard PC.

Fig. \ref{alldpol} concerns small values of dimension $d\le 4$.
It is seen that the solid lines form closed or semi-open curves which enclose 
disjunctive regions of the complex plane.
There are no critical zeros as $\rho_x=\frac{d}{2}\to 0$, in agreement with
the analysis of the previous section.
In the studied interval of $\rho_y$-values, there is always a finite gap between
the axis $\rho_x=0$ and the smallest $\rho_x$-component of the critical zeros;
an open question is whether the non-zero gap is present for all critical 
zeros (with larger values of $\rho_y$).

Fig. \ref{larged3} concerns large values of $d$, up to $d=50$,
where the asymptotic features of the critical zeros begin to occur.
The solid curve connecting the first zeros (with the smallest 
$\rho_y$-coordinate) ends up at $d=9.2455524667456$,
but the curve continues reflection-symmetrically across the $\rho_x$-axis
into the lower quadrant. 
This is why the critical zero 1b with the zero $\rho_y$-component
is a right edge point. 
The solid lines of the second and third zeros (as they appear in the region 
of small $d$) proceed up to an infinite dimension. 
For large values of $d$, the critical zeros are distributed 
equidistantly along the imaginary axis, as is predicted by the
asymptotic relation (\ref{infinity}) (the dashed lines).
The spacing between the nearest-neighbour zeros vanishes as $d\to\infty$.

From all critical zeros lying on a given curve, the ``edge'' points, 
denoted as $\rho^* = (d^*/2,\rho_y^*)$ with $\rho_y^*\equiv \rho_y(d^*)$, 
are the most relevant.
They are defined by the tangent
${\rm d}\rho_y/{\rm d}d\vert_{\rho^*} = \pm \infty$ 
or, equivalently, ${\rm d}d/{\rm d}\rho_y\vert_{\rho^*} = 0$.
They satisfy Eq. (\ref{standard}) for the critical zeros, i.e.,    
\begin{equation} \label{eq1}
- \frac{d^*}{\left(\frac{d^*}{2}\right)^2+{\rho_y^*}^2} + \int_0^1 {\rm d}t\, 
t^{d^*/4-1} \cos\left( \frac{\rho_y^* \ln t}{2} \right) \left[ 
\theta_3^{d^*}\left({\rm e}^{-\pi t}\right) - t^{-d^*/2} \right] = 0 , 
\end{equation} 
and simultaneously the derivative of Eq. (\ref{standard}) with respect to 
$\rho_y$, taken with ${\rm d}d/{\rm d}\rho_y\vert_{\rho^*} = 0$,
\begin{eqnarray} 
\frac{2 \rho_y^* d^*}{\left[ \left( \frac{d^*}{2}\right)^2+
{\rho_y^*}^2\right]^2} - \frac{1}{2} \int_0^1 {\rm d}t\, t^{d^*/4-1} 
\sin\left( \frac{\rho_y^* \ln t}{2} \right) \ln t  & & \nonumber \\ \times
\left[ \theta_3^{d^*}\left({\rm e}^{-\pi t}\right) - t^{-d^*/2} \right] 
& = & 0 . \label{eq2}
\end{eqnarray} 
The set of equations (\ref{eq1}) and (\ref{eq2}) has an infinite number
of solutions for $\rho^*$ with dimension $d^*$ being in general non-integer.
We have to distinguish between the ``left'' edge points, for which
${\rm d}^2d/{\rm d}\rho_y^2\vert_{\rho^*} > 0$, and the ``right'' edge points,
for which ${\rm d}^2d/{\rm d}\rho_y^2\vert_{\rho^*} < 0$.
The coordinates and left/right orientation of edge points in figures 
\ref{alldpol} and \ref{larged3} are summarized in Tab. \ref{table}.

\begin{table}[tbh]
\caption{The coordinates and the orientation of the edge points in Figs.
\ref{alldpol} and \ref{larged3}.}
\label{table}
\begin{tabular}{cccc}
\hline
edge point & orientation &$\rho^*_x=d^*/2$&$\rho^*_y$\\
\hline
1a &left&0.10846187908294 &18.06404476224324 \\
1b &right&4.62277623337280 & 0 \\
2a &left&0.029260757098957 &28.25989865119296 \\
2b &right&1.13615655471973 &27.06485479190591 \\
3a &left&0.076684964492103 &36.29956597219118 \\
3b &right&0.17608667918405 &38.97086173076263 \\
3c &left&0.023788974966443 &42.00296457563092 \\
3d &right&0.69958436750509 &42.29187347594789 \\
3e &left&0.28286847694364 &39.08036320922192 \\
4a &left&0.94484709689530 &42.43883096807280 \\
4b &right&1.87159485174678 &42.20920217767993 \\\hline
\end{tabular}
\vspace*{-4pt}
\end{table}

$\rho_y$ as the function of $d$ is singular (non-analytic) at the
edge points.
Let us consider, say, a left edge point with coordinates
$(d^*/2,\rho_y(d^*))$ and study an infinitesimal deviation from it 
along the curve of critical zeros $\rho_y(d)$.
Setting 
\begin{equation} \label{drho}
d=d^* + \Delta d \quad (0<\Delta d\ll 1), \quad
\rho_y(d) = \rho_y(d^*) + \Delta\rho_y \quad (\Delta\rho_y\ll 1)
\end{equation}
in (\ref{standard}) and expanding systematically in powers of small 
deviations $\Delta d$ and $\Delta\rho_y$, one gets
\begin{equation} \label{equ}
\alpha \Delta d + \gamma (\Delta\rho_y)^2 - \beta\Delta d \Delta\rho_y
- \delta (\Delta\rho_y)^3 
+ {\cal O}\left[(\Delta d)^2\right] +
{\cal O}\left[\Delta d(\Delta\rho_y)^2\right] = 0 , 
\end{equation} 
where the expansion coefficients are given by
\begin{eqnarray}
\alpha & = & \frac{\left(\frac{d^*}{2}\right)^2-{\rho_y^*}^2}{\left[
\left(\frac{d^*}{2}\right)^2+{\rho_y^*}^2\right]^2} +
\int_0^1 {\rm d}t\, t^{\frac{d^*}{4}-1} 
\cos\left( \frac{\rho_y^* \ln t}{2} \right) \nonumber \\ & & \times
\left\{ \frac{\ln t}{4} \left[ \theta^{d^*}_3\left({\rm e}^{-\pi t}\right) + 
t^{-d^*/2} \right] + \theta_3^{d^*}\left({\rm e}^{-\pi t}\right)
\ln\left[ \theta_3\left({\rm e}^{-\pi t}\right) \right] \right\} , 
\label{alpha} \\ 
\beta & = & 2 \rho_y^* \frac{3\left(\frac{d^*}{2}\right)^2-{\rho_y^*}^2}{
\left[\left(\frac{d^*}{2}\right)^2+{\rho_y^*}^2\right]^3} + \frac{1}{2}
\int_0^1 {\rm d}t\, t^{\frac{d^*}{4}-1} 
\sin\left( \frac{\rho_y^* \ln t}{2} \right) \ln t \nonumber \\ & & \times 
\left\{ \frac{\ln t}{4} \left[ \theta_3^{d^*}\left({\rm e}^{-\pi t}\right) + 
t^{-d^*/2} \right] + \theta_3^{d^*}\left({\rm e}^{-\pi t}\right)
\ln\left[ \theta_3\left({\rm e}^{-\pi t}\right) \right] \right\} , 
\label{beta} \\ 
\gamma & = & d^* \frac{\left(\frac{d^*}{2}\right)^2-3{\rho_y^*}^2}{
\left[\left(\frac{d^*}{2}\right)^2+{\rho_y^*}^2\right]^3} 
- \frac{1}{8} \int_0^1 {\rm d}t\, t^{\frac{d^*}{4}-1} 
\cos\left( \frac{\rho_y^* \ln t}{2} \right) (\ln t)^2 \nonumber \\ & & \times 
\left[ \theta_3^{d^*}\left({\rm e}^{-\pi t}\right) - t^{-d^*/2} \right] , 
\label{gamma} \\ 
\delta & = & 4 d^* \rho_y^* \frac{\left(\frac{d^*}{2}\right)^2-{\rho_y^*}^2}{
\left[\left(\frac{d^*}{2}\right)^2+{\rho_y^*}^2\right]^4} - \frac{1}{48}
\int_0^1 {\rm d}t\, t^{\frac{d^*}{4}-1} 
\sin\left( \frac{\rho_y^* \ln t}{2} \right) (\ln t)^3 \nonumber \\ & & \times 
\left[ \theta_3^{d^*}\left({\rm e}^{-\pi t}\right) - t^{-d^*/2} \right] . 
\label{delta}
\end{eqnarray}
Note that the linear term of order $\Delta\rho_y$ is missing in (\ref{equ})
due to the validity of Eq. (\ref{eq2}) for the critical edge zeros.

For $d>d^*$ $(\Delta d>0)$, Eq. (\ref{equ}) tells us that the leading
order of the expansion of $\Delta\rho_y$ with respect to $\Delta d$ is given by
the relation $\alpha \Delta d + \gamma (\Delta\rho_y)^2 = 0$, i.e.,
\begin{equation} \label{singular}
\Delta \rho_y \sim \pm \sqrt{-\frac{\alpha}{\gamma}} \sqrt{\Delta d} ,
\end{equation}
where the prefactor sign $\pm$ specifies the up and down branches
of the plot $\rho_y(d)$.
Notice that $-\alpha/\gamma$ must be a positive number to get a real solution
for $\rho_y(d)$ and our numerical calculations confirm that for all studied 
left edge points it really is so.
Singular relations of type (\ref{singular}) with exponent $\frac{1}{2}$ occur 
also for the order parameter in critical phenomena of statistical systems at 
the second-order phase transition within the so-called mean-field approach 
\cite{Baxter82,Samaj13}.
The next order of the expansion of $\Delta\rho_y$ in $\Delta d$ follows from 
Eq. (\ref{equ}) by considering two additional terms of the order 
$(\Delta d)^{3/2}$.
Adding a term $c \Delta d$ to the formula for $\Delta \rho_y$
in (\ref{singular}) and expanding all functions up to the order
$(\Delta d)^{3/2}$ fixes $c$ as follows 
$2\gamma c = \beta - \alpha \delta/\gamma$. 
We conclude that
\begin{eqnarray} 
\rho_y(d) & = & \rho_y(d^*) \pm \sqrt{-\frac{\alpha}{\gamma}} \sqrt{d-d^*}
+ \frac{1}{2\gamma} \left( \beta - \frac{\alpha\delta}{\gamma} \right)
(d-d^*) \nonumber \\ & &
+ {\cal O}\left[(d-d^*)^{3/2}\right] . \label{rhoypositive}
\end{eqnarray}
A similar analysis can be made for right edge zeros.

\section{Generation of off-critical zeros from critical edge zeros} \label{Sec6}
This section is about a continuous generation of off-critical zeros 
from the critical edge zeros. 

In the case of the left edge zero, Eq. (\ref{equ}) with $-\alpha/\gamma>0$
has no real solution for $\Delta\rho_y$ if $\Delta d = d-d^* <0$.
Let us assume that for $\Delta d<0$ there is also a continuous deviation of 
the $\rho_x$-component from its critical value $\frac{d}{2}$, see Eq.
(\ref{Deltarho}).
Considering then the relations (\ref{drho}) in equations (\ref{nonstandard1}) 
and (\ref{nonstandard2}) and expanding in powers of small variables 
$\Delta d$, $\Delta\rho_x$ and $\Delta\rho_y$, one obtains
\begin{eqnarray} 
\alpha \Delta d - \gamma (\Delta\rho_x)^2 + \gamma (\Delta\rho_y)^2 
- \beta\Delta d \Delta\rho_y + 3\delta (\Delta\rho_x)^2 \Delta\rho_y & &
\nonumber \\ 
- \delta (\Delta\rho_y)^3 + \cdots & = & 0 , \label{firsty} \\
\Delta\rho_x \left[ \beta \Delta d - 2\gamma \Delta\rho_y
- \delta (\Delta\rho_x)^2 + 3\delta (\Delta\rho_y)^2
+ \cdots \right] & = & 0 . \label{secondy}
\end{eqnarray} 

As is evident from the first relation (\ref{firsty}), setting $\Delta\rho_x=0$, 
$\Delta\rho_y$, which is real for $\Delta d>0$, becomes pure imaginary for 
$\Delta d<0$ (in the leading order of its expansion in $\Delta d$)
which is in contradiction with the definition of $\rho_y$ as a real number.
However, the term $\gamma (\Delta\rho_y)^2$ containing the variable 
$\Delta\rho_y$ has a counterpart with the opposite sign 
$-\gamma (\Delta\rho_x)^2$ containing the variable $\Delta\rho_x$ 
as the actual candidate for the leading-order symmetry breaking 
$\Delta\rho_x\ne 0$. 
In the leading order, the equation 
$\alpha \Delta d - \gamma (\Delta\rho_x)^2 = 0$ implies that 
\begin{equation} \label{singular1}
\Delta \rho_x \sim \pm \sqrt{-\frac{\alpha}{\gamma}} \sqrt{-\Delta d} ,
\end{equation}
where the $\pm$ sign reflects the split of $\Delta\rho_x$ onto
two different branches (see discussion below).  
Since $\Delta\rho_x\ne 0$, the second relation (\ref{secondy}) implies 
that the leading order is determined for $\Delta\rho_y$ by
$\beta \Delta d - 2\gamma\Delta\rho_y -\delta (\Delta\rho_x)^2=0$, i.e., 
\begin{equation} \label{rhoynegative}
\Delta\rho_y = - \frac{1}{2\gamma} 
\left( \beta - \frac{\alpha\delta}{\gamma} \right) (d^*-d) 
+ {\cal O}\left[(d^*-d)^{3/2}\right] , \qquad d<d^*.
\end{equation}
The next order of the expansion of $\Delta\rho_x$ in $-\Delta d$ follows from 
Eq. (\ref{firsty}) by adding a term $c (-\Delta d)$ to the formula for
$\Delta \rho_x$ in (\ref{singular1}).
The term arising from $(\Delta\rho_x)^2$, which is proportional to
$c(-\Delta d)^{3/2}$, has no counterparts since all additional terms in
(\ref{firsty}) are of the order $(-\Delta d)^2$.
Consequently, $c=0$ and one gets
\begin{equation} \label{rhoxnegative}
\Delta\rho_x = \pm \sqrt{-\frac{\alpha}{\gamma}} \sqrt{d^*-d}
+ {\cal O}\left[(d^*-d)^{3/2}\right] , \qquad d<d^*.
\end{equation}
We conclude that the expansion of up and down branches of $\rho_y(d)$ 
(\ref{rhoypositive}), valid for $d>d^*$, splits into the expansion for 
the left and right branches (tails) (\ref{rhoxnegative}) and 
(\ref{rhoynegative}), valid for $d<d^*$, with similar expansion 
coefficients.
While the previous two up and down branches of $\rho_y(d)$ for $d>d^*$
were not bound by a symmetry, the positive and negative branches
of off-critical zeros $[\rho_x(d),\rho_y(d)]$ for $d<d^*$ are bound by 
the symmetry $[\rho_x(d),\rho_y(d)]\to [d-\rho_x(d),\rho_y(d)]$ which 
takes place for every deviation from the edge point.
A similar analysis can be made for right edge zeros, to keep the parameter 
$\Delta d$ negative from the side of off-critical zeros it must be defined as 
$\Delta d = d^*-d$.

The numerical evaluation of off-critical zeros using \textit{Mathematica} 
is relatively simple due to the continuity of tails as $d$ starts
to deviate from $d^*$. 
As the function to deal with, one takes the sum of two squares of the lhs 
of equations (\ref{nonstandard1}) and (\ref{nonstandard2}) which 
should vanish at zeros.
The command FindMinimum is applied to the function and the zero is taken as 
guaranteed if the function value is of order at most $10^{-23}$.
To prevent escape from the local minimum, one starts from 
(say right) edge point by shifting the dimension by a tiny amount 0.0001,
after few steps the shift can be increased to 0.001 and finally to 0.01.
To find a minimum takes approximately 60 seconds of CPU time on a standard PC.

In what follows, various scenarios are presented how the two conjugate tails of 
off-critical zeros behave during their $d$-evolution. 

\begin{figure}[th]
\centering
\includegraphics[clip,width=0.75\textwidth]{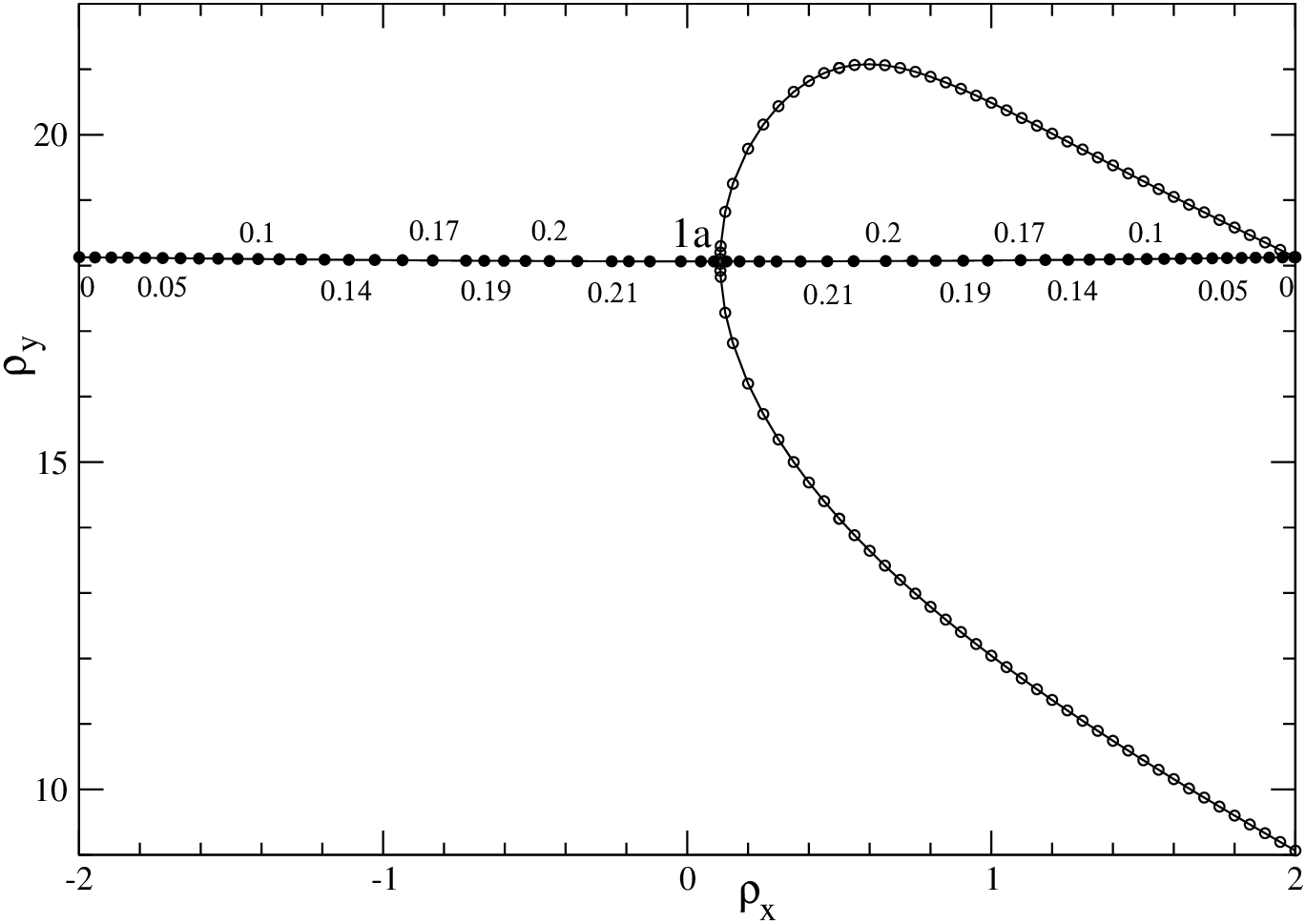}
\caption{The left and right tails of off-critical zeros (full circles)
starting from the critical edge point 1a $(d^*\approx 0.217)$.
The dimension of off-critical zeros is indicated at a few points.
The tails go down to $d=0$ at the points $(\pm 2,4\pi/\ln 2)$.}
\label{chvost1}
\end{figure}

The simplest scenario, presented in Fig. \ref{chvost1}, is associated with 
the left edge point denoted as 1a in Fig. \ref{alldpol}. 
The critical zeros are represented by open circles, the left and right
tails of off-critical zeros by full circles.
The dimension of off-critical zeros is indicated at a few points.
The left and right tails of off-critical zeros start from $d^*\approx 0.217$
at the edge point 1a and go down monotonously to $d=0$ at points 
$(-2,4\pi/\ln 2)$ and $(2,4\pi/\ln 2)$, respectively.  
The latter points coincide with the $k=1$ zeros (\ref{sk}) found in $d=0$.
As is seen in the figure, the end-point of the right tail $(2,4\pi/\ln 2)$
coincides with a critical zero at $d=4$.

\begin{figure}[th]
\centering
\includegraphics[clip,width=0.75\textwidth]{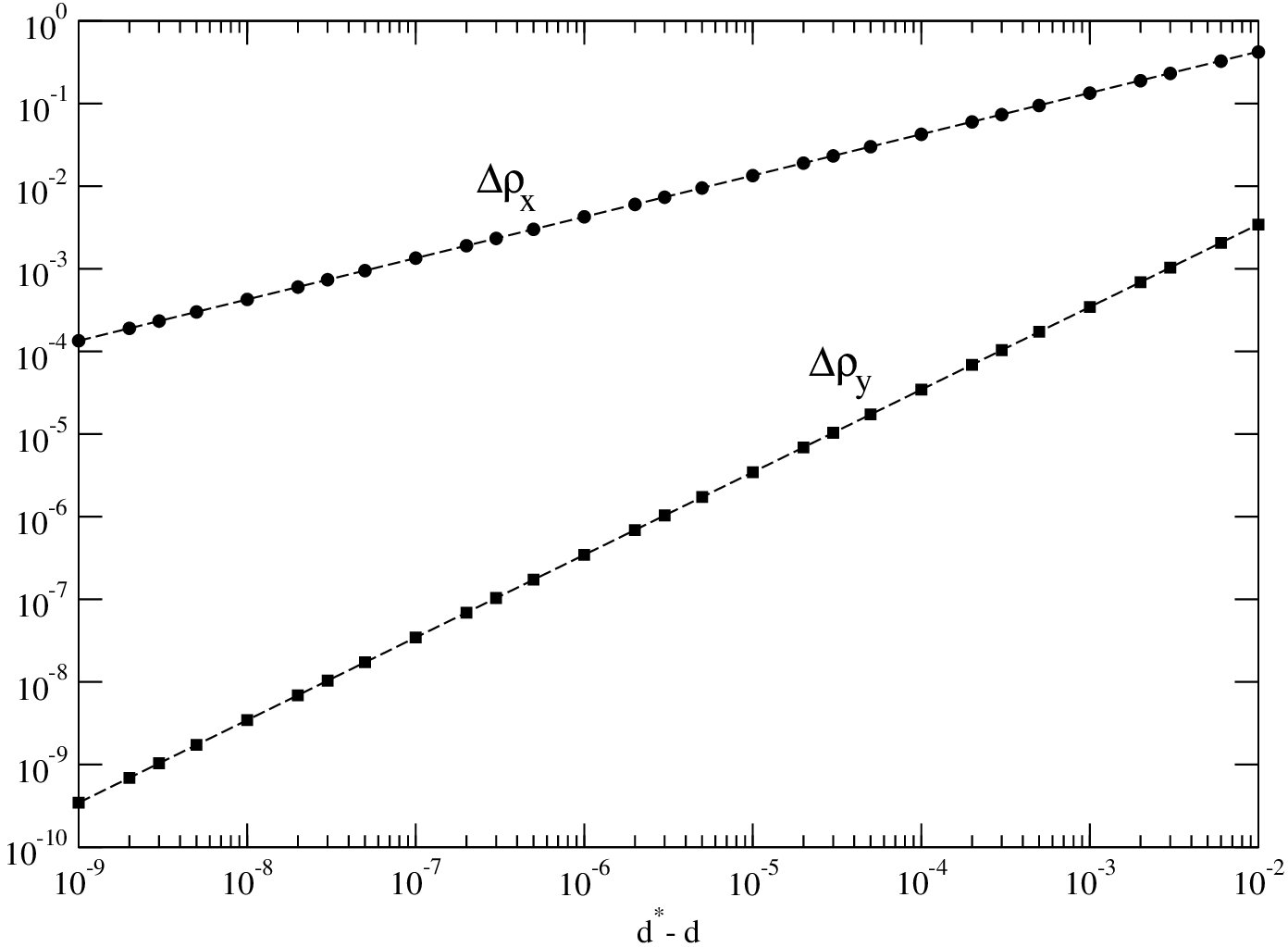}
\caption{The log-log plot of numerical data for the dependence of 
$\Delta\rho_x$ (full circles) and $\Delta\rho_y$ (full squares)
on the small dimension deviation $d^*-d$ for the right tail 
going from the left edge point 1a.
Data fit perfectly with the analytic predictions (\ref{edge1a})
represented in logarithmic scale by the dashed straight lines.}
\label{crit}
\end{figure}

The test of the expansion formulas (\ref{rhoynegative}) and 
(\ref{rhoxnegative}) for the right tail generated from 
the left edge point 1a is presented in Fig. \ref{crit}. 
Our numerical data for the dependence of $\Delta\rho_x$ and $\Delta\rho_y$ 
on $d^*-d$ are represented in logarithmic scale by full circles 
and squares, respectively. 
For small deviations from the edge point $d^*-d\to 0^+$, the expansion 
formulas (\ref{rhoynegative}) and (\ref{rhoxnegative}), with the
constants $\alpha,\beta,\gamma,\delta$ evaluated by using the formulas
(\ref{alpha})-(\ref{delta}), lead to
\begin{equation} \label{edge1a}
\Delta\rho_x \sim \pm 4.24563 \sqrt{d^*-d} , \qquad 
\Delta\rho_y \sim 0.344516 (d^*-d) .
\end{equation}
The log-log plots of these analytic predictions, represented in 
Fig. \ref{crit} by the dashed straight lines, fit perfectly the corresponding 
numerical data for the small deviation $d^*-d$ ranging 
from $10^{-9}$ to $10^{-2}$.  

\begin{figure}[th]
\centering
\includegraphics[clip,width=0.75\textwidth]{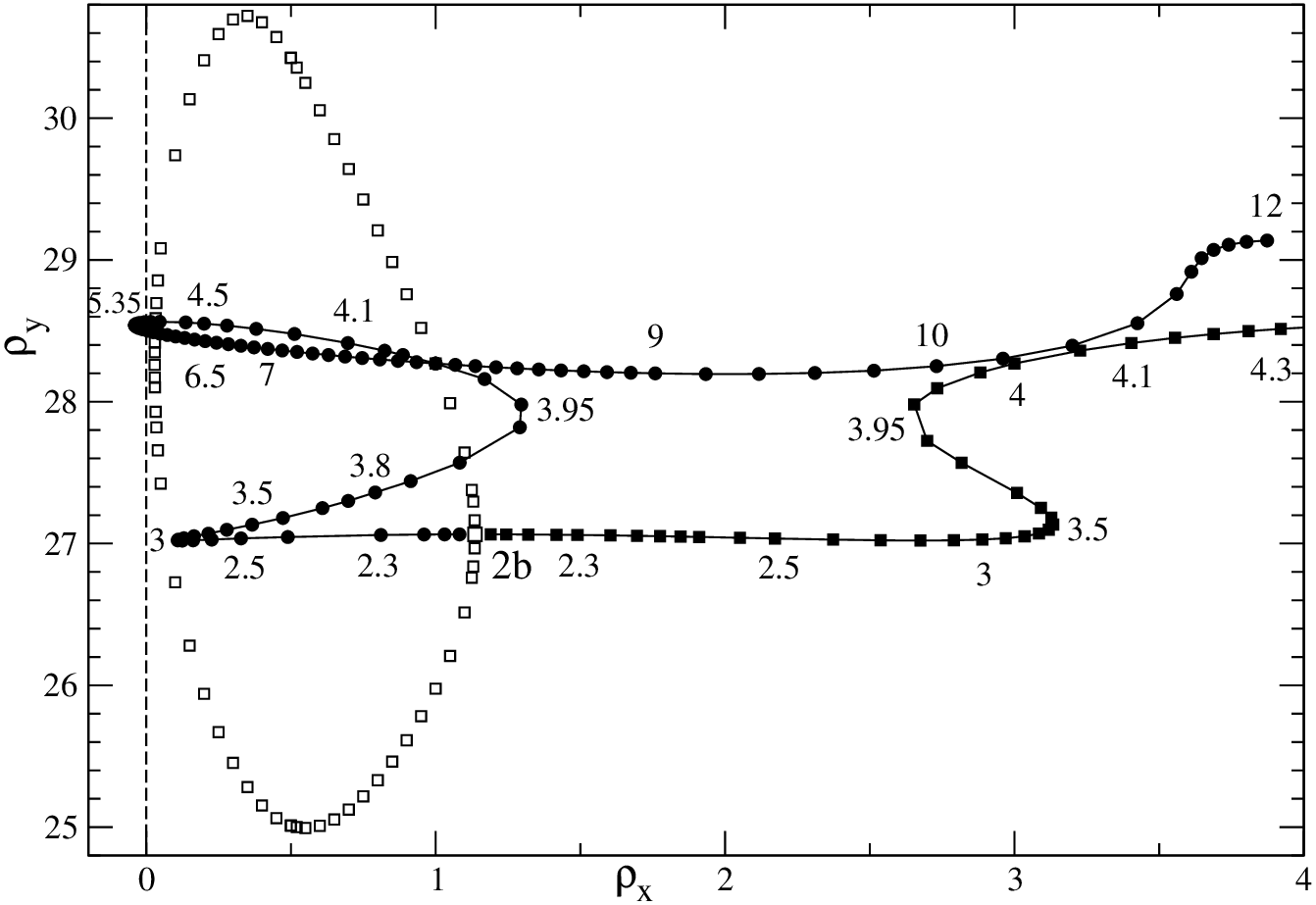}
\caption{The left (full circles) and right (full squares) tails of off-critical
zeros starting from the critical edge point 2b $(d^*\approx 2.272)$.
Dimension $d$ increases along the tails up to $\infty$, only fragments 
of the tails are shown in the figure.} 
\label{zeroschvosty}
\end{figure}

\begin{figure}[th]
\centering
\includegraphics[clip,width=0.75\textwidth]{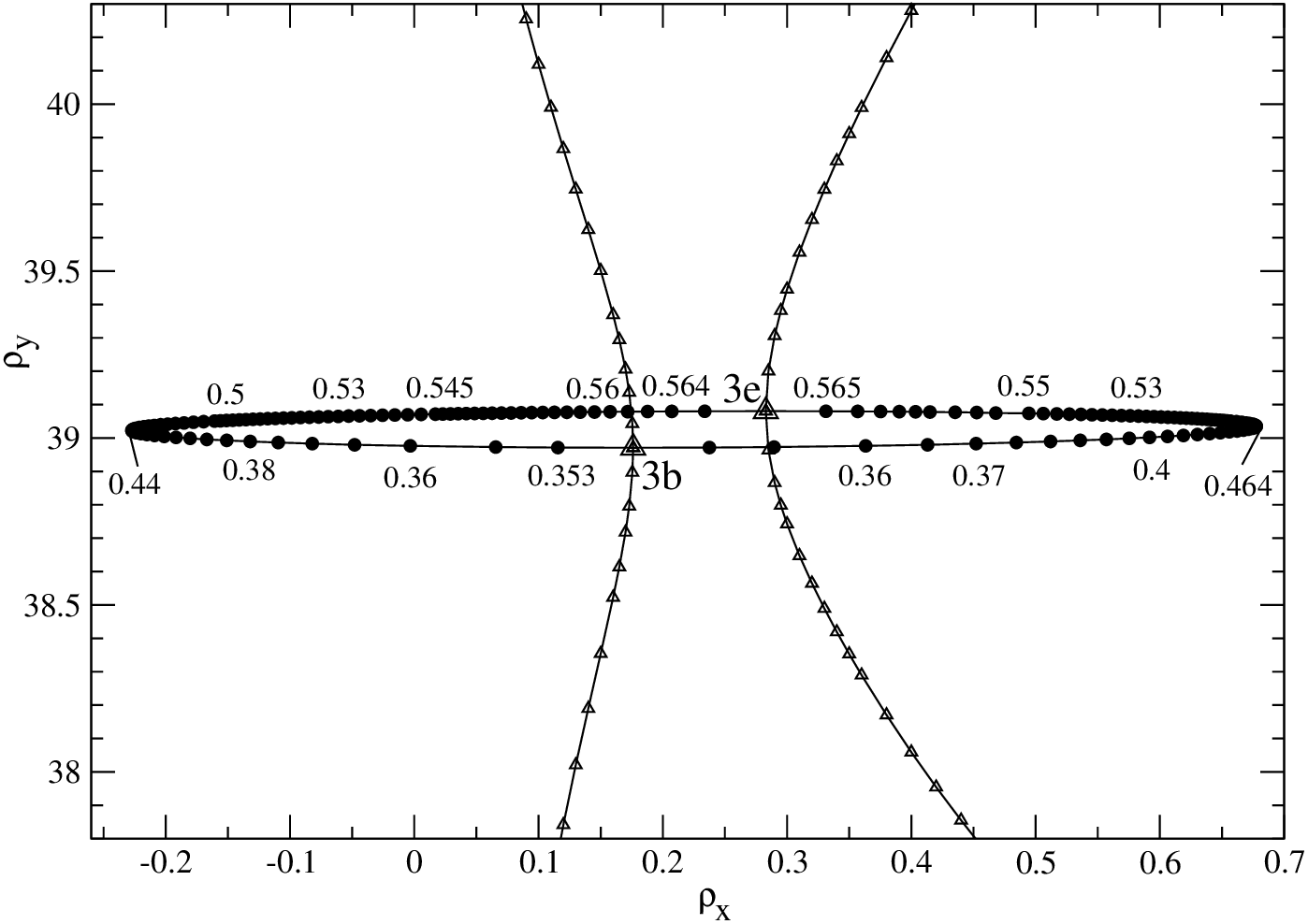}
\caption{Interpolation of off-critical zeros (full circles) between 
the right edge point 3b $(d^*\approx 0.352)$ and the left edge point 3e 
$(d^*\approx 0.566)$, both edge points lying on the same curve 
of critical zeros (open triangles).}
\label{hodinypop}
\end{figure}

The form of the tails of off-critical zeros is more complicated in the case of 
the right edge point 2b (as denoted in Fig. \ref{alldpol}), see
Fig. \ref{zeroschvosty}.
Because of the right orientation of this edge point, dimension $d$ increases 
along the left (full circles) and right (full squares) tails from 
$d^*\approx 2.272$ to $\infty$. 
As a check, for the left tail we recover a special off-critical zero 
which occurs simultaneously in dimensions $d=4$ and $d=8$, see 
equations (\ref{zeta4}) and (\ref{zeta8}), and coincide with 
the fourth critical zero of the $d=2$ Epstein zeta-function.
For the dimension of special physical interest $d=3$, the left tail contains 
the off-critical zero
$\rho\approx 0.111189793551259+27.0278811412527548{\rm i}$.
It is interesting that the left tail covers also a tiny interval of negative
values of $\rho_x$ for dimensions from the interval $[4.965,5.775]$.
This means that for integer dimension $d=5$ there exists a pair of conjugate 
off-critical zeros with the same imaginary component, namely 
$\rho_y\approx 28.5599914110240345$, and the real components 
$\rho_x\approx -0.00717997528701, 5.00717997528701$ which lie outside 
the critical strip $[0,5]$.
This phenomenon is usually present in lower dimensions with critical
strips of small width, especially for $d$ between 0 and 1,
see Figs. \ref{chvost1} and \ref{hodinypop}. 

Another scenario is presented in Fig. \ref{hodinypop} where the two 
tails of off-critical zeros (full circles) interpolate between 
the right edge point 3b $(d^*\approx 0.352)$ and the left edge point 3e 
$(d^*\approx 0.566)$, both points lying on the same curve of critical zeros
(open triangles).  
Although the interval of dimensions is relatively narrow, the interval of 
$\rho_x$-values of off-critical zeros is relatively large and involves also 
negative values.

\begin{figure}[th]
\centering
\includegraphics[clip,width=0.75\textwidth]{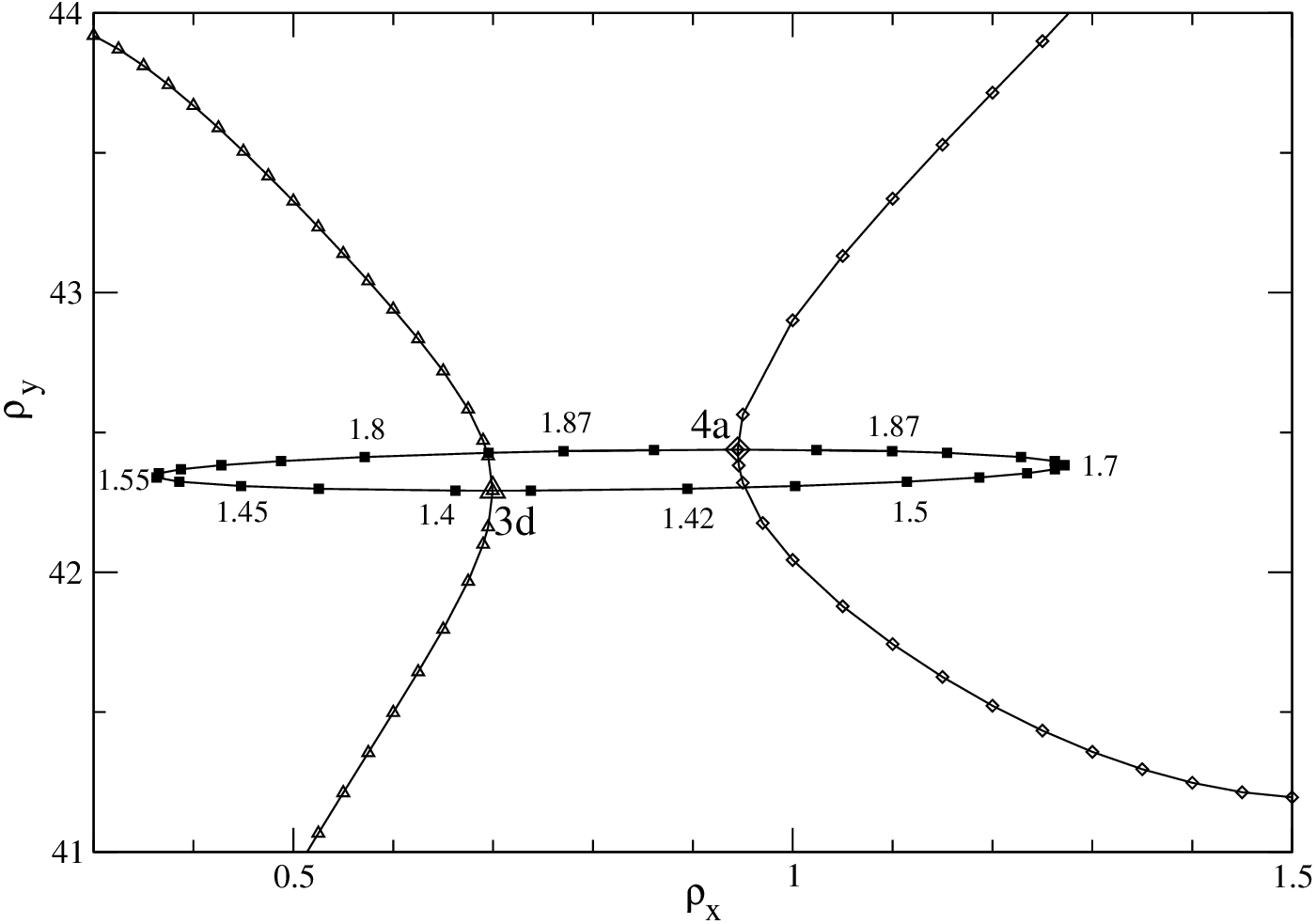}
\caption{Interpolation of off-critical zeros (full squares) between 
the right edge point 3d $(d^*\approx 1.399)$ and the left edge point 4a 
$(d^*\approx 1.890)$, the edge points lying on distinct curves 
of critical zeros (open triangles and diamonds).}
\label{horslza}
\end{figure}

In Fig. \ref{horslza}, the two tails of off-critical zeros 
(full squares) interpolate between the right edge point 3d $(d^*\approx 1.399)$ 
and the left edge point 4a $(d^*\approx 1.890)$ which lie on different 
curves of critical zeros (open triangles and diamonds).  

\begin{figure}[th]
\centering
\setbox1=\hbox{\includegraphics[clip,height=7.8cm]{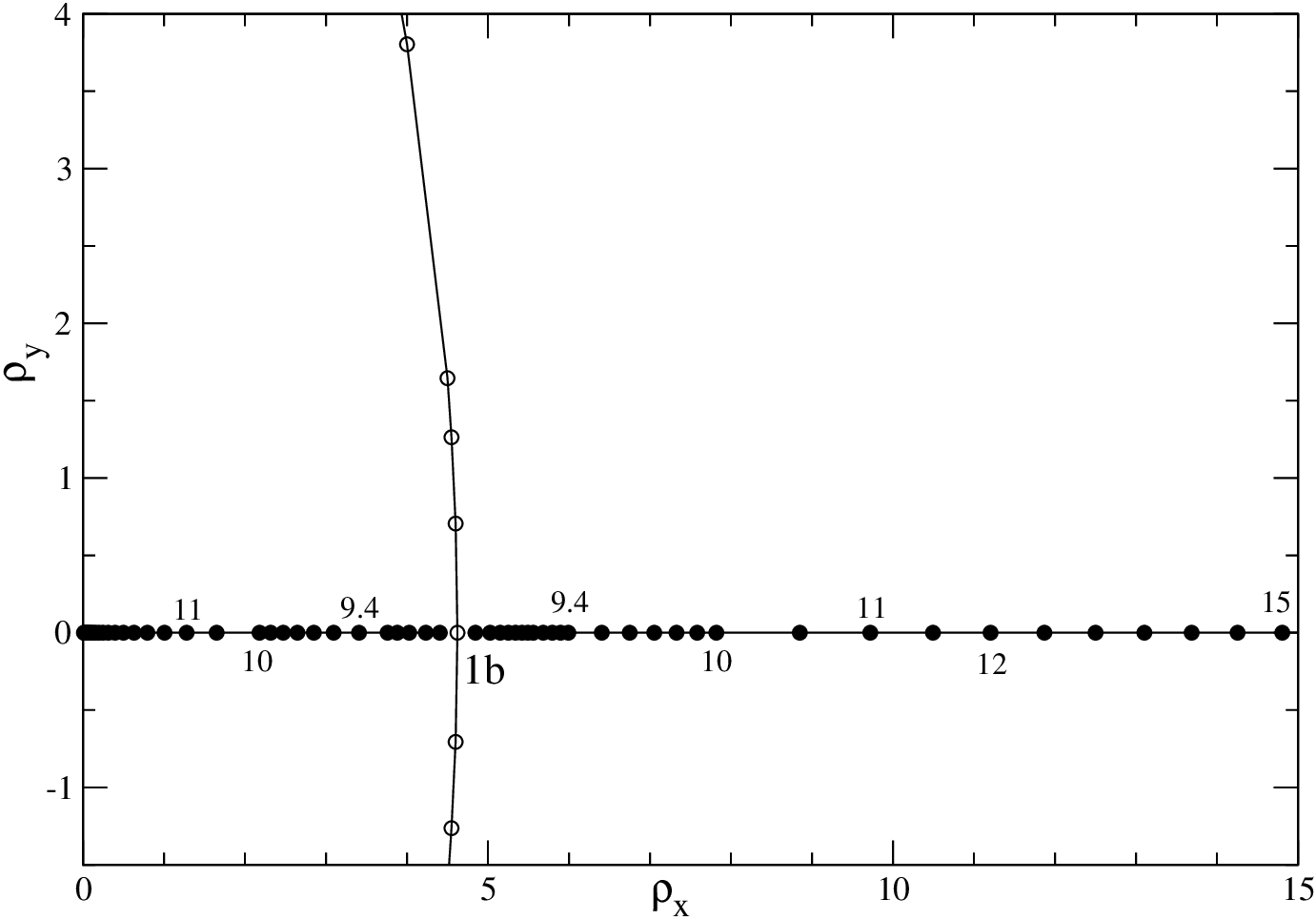}}
\includegraphics[clip,height=7.8cm]{fig8a.eps}\llap{
\raisebox{3.59cm}{\includegraphics[clip,height=4.1cm]{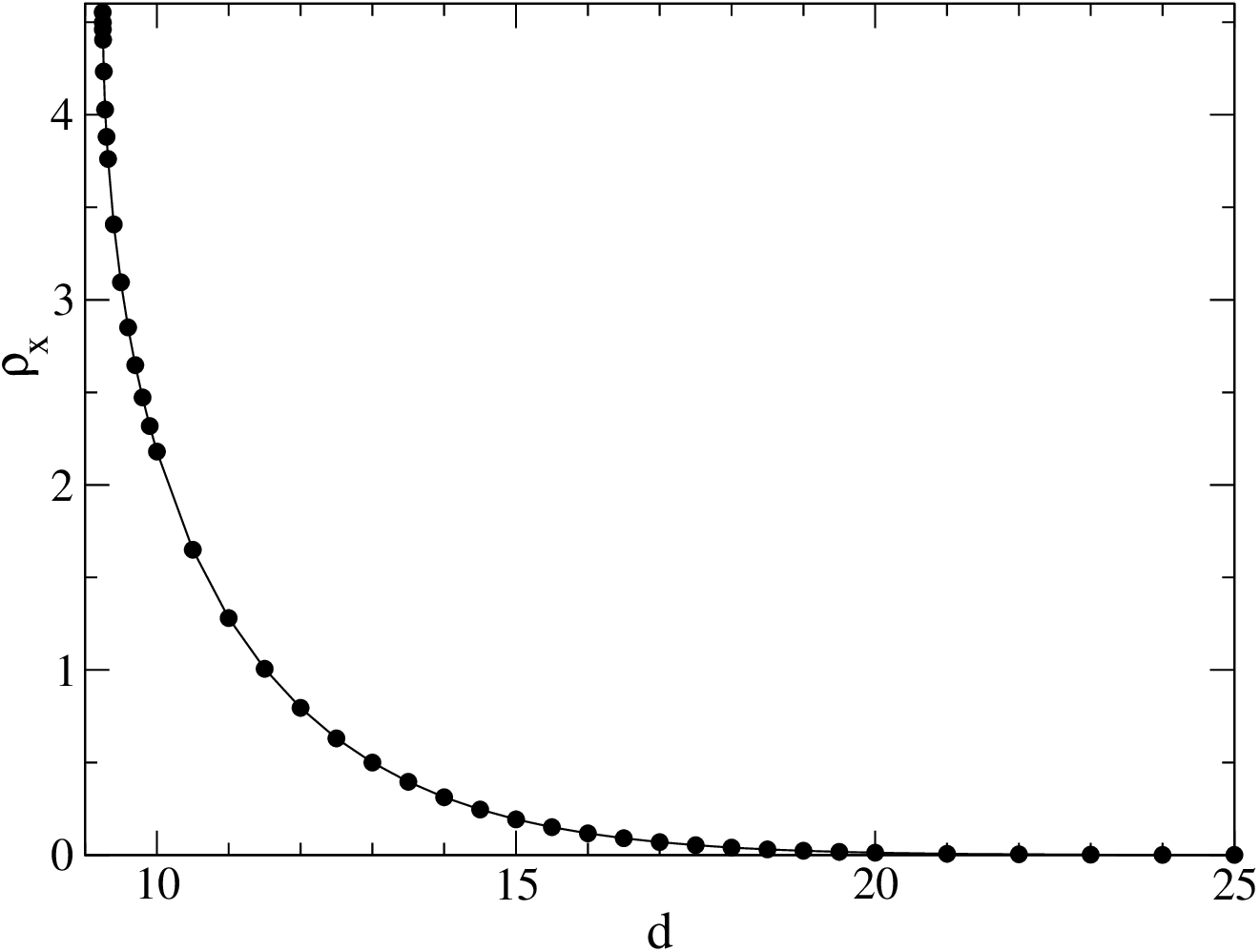}}}
\caption{Two conjugate tails of the real off-critical zeros (full circles) 
coming from the right edge point 1b $(d^* \approx 9.24555)$, 
localized at the intersection of the curve of critical zeros (open circles) 
and the real axis.
The dimension of off-critical zeros is indicated at a few points.
The inset describes the quick approach of the numerical values of
$\rho_x$ (full circles) to zero when increasing dimension $d$.}
\label{ep1b}
\end{figure}

The critical point 1b in Fig. \ref{larged3} is a right edge point
because the curve continues reflection-symmetrically across
the $\rho_x$-axis into the lower quadrant. 
The position of this point $({\rho_c^*})_x = d_c^*/2$ can be found 
by setting $\rho_y=0$ in (\ref{standard}),
\begin{equation}
- \frac{4}{d_c^*} + \int_0^1 {\rm d}t\, t^{d_c^*/4-1} \left[ 
\theta_3^{d_c^*}\left({\rm e}^{-\pi t}\right) - t^{-d^*_c/2} \right] = 0 . 
\end{equation}  
The numerical solution of this equation is $d_c^* = 9.24555\ldots$.
For every $d>d_c^*$, there exists a pair of off-critical zeros on 
the real axis [note that Eq. (\ref{nonstandard2}) holds automatically 
for $\rho_y=0$] whose conjugate components $\rho_x$ and $d-\rho_x$ satisfy 
the integral equation
\begin{equation} \label{expr}
- \frac{1}{\rho_x} - \frac{1}{d-\rho_x} + \int_0^1 \frac{{\rm d}t}{t}
\frac{t^{\rho_x/2}+t^{(d-\rho_x)/2}}{2}
\left[ \theta_3^d\left({\rm e}^{-\pi t}\right) - t^{-d/2} \right] = 0 .   
\end{equation}
The two solutions of this equation are lying in the critical strip 
$0<\rho_x<d$, see Fig. \ref{ep1b}. 
This can be explained for integral dimensions $d$ by the fact that for real 
$\rho>d$ the Epstein zeta-function is the sum of positive numbers and 
therefore cannot vanish.
To be more particular, for $d=10$ one has the pair of real zeros
$\rho=2.17985543147...;7.82014456853...$, for $d=12$ one has 
$\rho=0.7951625733...;11.2048374267...$, for $d=20$ one has 
$\rho=0.0127182144...;19.9872817856...$, etc.
In the large-$d$ limit, the two solutions tend to the boundaries 
$0$ and $d$ of the critical strip.
The quick approach to 0 of the numerical values of $\rho_x$ (full circles)
is presented in the inset of Fig. \ref{ep1b}.

\section{Conclusion} \label{Sec7}
The basic definition of the hypercubic Epstein zeta-function 
(\ref{zetad}) requires an integer value of the spatial dimension $d$.
The analytic continuation of the lattice sum to the whole complex $s$-plane 
(\ref{zetadbest}) is well defined also for non-integer values of $d$.  
This extension was used to obtain numerically the closed or semi-open curves 
$\rho_y(d)$ of critical zeros (on the critical line), see figures 
\ref{alldpol} and \ref{larged3}. 
Each curve involves a finite number of left/right critical edge points 
$\rho^*=(d^*/2,\rho_y(d^*))$, defined by an infinite tangent 
${\rm d}\rho_y/{\rm d}d\vert_{\rho^*}$.
The coordinate data in Tab. \ref{table} indicate that the dimension of
edge critical points is in general non-integer.
As was shown in section \ref{Sec5}, the function $\rho_y(d)$ is a singular
function of the dimension deviation $d-d^*>0$ $(d^*-d>0)$ for the left
(right) edge points.
Changing the sign of the dimension deviation to $d-d^*<0$ $(d^*-d<0)$ 
for the left (right) edge points, these edge points give rise to two 
conjugate tails of off-critical zeros (off the critical line) with 
continuously varying dimension $d$, in formal analogy with critical 
phenomena for many-body statistical systems (section \ref{Sec6}).
Various versions of the generation mechanism are presented. 
Fig. \ref{chvost1} documents the generation of the left and right tails
of off-critical zeros from the left edge point 1a $(d^*\approx 0.217)$,
the dimension along the tails goes down to $0$ at the off-critical zeros 
$(\pm 2,4\pi/\ln 2)$.
The corresponding log-log plots of numerical data for
$\Delta\rho_x(d) = \rho_x(d)-\frac{d^*}{2}$ and 
$\Delta\rho_y(d) = \rho_y(d)-\rho_y(d^*)$ in the case of the right tail,
presented in Fig. \ref{crit}, are in perfect agreement with the analytic
prediction (\ref{edge1a}) valid for small dimension deviations $d^*-d$.
The generation of the off-critical tails from the right edge point 2b,
with dimension along tails going up to infinity, is pictured in
Fig. \ref{zeroschvosty}.
The interpolation of off-critical zeros between the right edge point 3b
and the left edge point 3e, both edge points lying on the same curve
of critical zeros, is presented in Fig. \ref{hodinypop}.
Fig. \ref{horslza} concerns an interpolation of off-critical zeros between
the right edge point 3d and the left edge point 4a, the edge points lying
on different curves of critical zeros.
For every $d>d^*_c\approx 9.246$, there exists a pair of conjugate
off-critical zeros on the real axis, having their origin in the right edge 
point 1b, see Fig. \ref{ep1b}.
As $d\to\infty$, the two zeros tend very quickly to the boundaries $0$ and $d$ 
of the critical strip.

As a by-product of the formalism, we have derived the exact formula (\ref{fr})
for $\lim_{d\to 0^+}\zeta^{(d)}(s)/d$.
This formula tells us that there are no critical zeros in the limit $d\to 0^+$.
In the studied interval of $\rho_y$-values smaller than 45, there is always 
a finite gap between the axis $\rho_x=0$ and the smallest $\rho_x$-component of 
the critical zeros; a remaining open question is whether a non-zero gap 
is present for all critical zeros (with larger values of $\rho_y$).
The spectrum of off-critical zeros in the limit $d\to 0^+$ was checked to
fulfill the obligatory sum rules.
Another check of the spectrum is that off-critical tails generated
from the left edge points end correctly at the $d\to 0^+$ off-critical zeros.
The exact treatment of the large-$d$ limit in section \ref{deinfinity}   
predicts an equidistant distribution of the critical zeros along the imaginary
axis (\ref{infinity}).
This result is confirmed numerically in Fig. \ref{larged3} where
the critical zeros (open circles, triangles and squares) approach to 
for sufficiently large $d$ the equidistant distribution (\ref{infinity})
represented by the dashed lines.

Another open question is whether the presented mechanism of generation of 
the tails of off-critical zeros from the critical edge points is the only one.
We anticipate that it is so.

\section*{Acknowledgements}
The support received from VEGA Grant No. 2/0092/21 
and Project EXSES APVV-16-0186 is acknowledged.

\end{document}